\documentclass[onefignum,onetabnum]{siamart220329}

\headers{FRAMEWORK FOR IMPLEMENTING VIRTUAL ELEMENT SPACES}{A. Dedner and A. Hodson}

\title{A framework for implementing general virtual element spaces
}

\author{Andreas Dedner\thanks{Department of Mathematics, University of Warwick, Coventry, CV4 7AL, UK 
  (\email{a.s.dedner@warwick.ac.uk}).}
\and Alice Hodson\thanks{Department of Numerical Mathematics, Charles University, Sokolovsk\'{a} 83, 186 75 Praha 8, Czech Republic 
(\email{hodson@karlin.mff.cuni.cz}).}}

\usepackage{amsmath}
\usepackage[amsmath]{ntheorem}
\usepackage{amssymb}
\usepackage{enumerate}
\usepackage{enumitem}
\usepackage{mathtools}
\usepackage[nodayofweek]{datetime} 
\usepackage[framemethod=tikz]{mdframed}
\usepackage{tikz}
\usepackage{tikz-cd}
\usepackage{listings}
\usepackage[draft]{todonotes} 
\usepackage{comment}
\usepackage{booktabs}
\usepackage{float}
\usepackage{multirow}
\usepackage{placeins}
\usepackage{ifthen}
\usepackage{pgfplots}
\usepackage{siunitx}
\usepackage{thmtools}
\usepackage[caption=false]{subfig}
\usepackage{algorithmic}
\usepackage[T1]{fontenc}
\usepackage{xcolor}
\usepackage[normalem]{ulem}
\usepackage{graphicx}

\pgfplotsset{compat=newest} 
\usepgfplotslibrary{units} 

\allowdisplaybreaks 

\numberwithin{equation}{section} 

\newsiamremark{assumption}{Assumption}
\newsiamremark{remark}{Remark}
\newsiamremark{example}{Example}

\newsiamremark{continuancex}{Continued Example}

\newenvironment{continuance}[1]
  {\edef\continuanceref{\getrefnumber{#1}}\begin{continuancex}}
  {\end{continuancex}}

\definecolor{codegreen}{rgb}{0,0.6,0}
\definecolor{codegray}{rgb}{0.5,0.5,0.5}
\definecolor{codepurple}{rgb}{0.58,0,0.82}
\definecolor{backcolour}{rgb}{0.95,0.95,0.92}

\lstdefinestyle{mystyle}{
    backgroundcolor=\color{backcolour},   
    commentstyle=\color{codegreen},
    keywordstyle=\color{magenta},
    numberstyle=\tiny\color{codegray},
    stringstyle=\color{codepurple},
    basicstyle=\ttfamily\footnotesize,
    breakatwhitespace=false,         
    breaklines=true,                 
    captionpos=b,                    
    keepspaces=true,                 
    numbers=left,                    
    numbersep=5pt,                  
    showspaces=false,                
    showstringspaces=false,
    showtabs=false,                  
    tabsize=2
}

\newcommand{\half}{\frac{1}{2}}
\newcommand{\R}{\mathbb{R}}

\newcommand{\N}{\mathbb{N}}

\newcommand{\prob}{\mathbb{P}}

\newcommand{\cB}{\mathcal{B}}

\newcommand{\cF}{\mathcal{F}}

\newcommand{\cM}{\mathcal{M}}
\newcommand{\cC}{\mathcal{C}}

\newcommand{\cP}{\mathcal{P}}
\newcommand{\cS}{\mathcal{S}}
\newcommand{\cT}{\mathcal{T}}

\newcommand{\cV}{\mathcal{V}}

\newcommand{\element}{E}

\newcommand{\polOrder}{\ell}

\newcommand{\valueProj}{\Pi^{\element}_0}
\newcommand{\gradProj}{\Pi^{\element}_1}
\newcommand{\hessProj}{\Pi^{\element}_2}
\newcommand{\edgeProj}{\Pi^s_0}
\newcommand{\faceNormalProj}{\Pi^s_n}
\newcommand{\faceProj}{\Pi^s_0}
\newcommand{\meshProj}{\Pi^F_0}

\newcommand{\vertiii}[1]{{\left\vert\kern-0.25ex\left\vert\kern-0.25ex\left\vert #1 
    \right\vert\kern-0.25ex\right\vert\kern-0.25ex\right\vert}}

\newcommand{\mesh}{\cT_h}

\newcommand{\dx}{\mathrm{d}x}
\newcommand{\ds}{\mathrm{d}s}

\newcommand{\valueBasisSet}{\cB_0}
\newcommand{\gradBasisSet}{\cB_1}
\newcommand{\hessBasisSet}{\cB_2}

\newcommand{\elementDofSet}{\Lambda^{\element}}
\newcommand{\interfaceDofSet}{\Lambda^{s}}

\newcommand{\edgeSet}{\cS_h(\element)}
\newcommand{\NoLocalDofs}{N_{dof}^{\element}}

\newboolean{thesis}
\newboolean{thesiscorrections}
\newboolean{changes}
\newboolean{examplechanges}
\newboolean{refereecommentsvisible}

\begin{document}

\maketitle

\setboolean{thesis}{false}
\setboolean{thesiscorrections}{false}
\setboolean{changes}{false}
\setboolean{examplechanges}{false}
\setboolean{refereecommentsvisible}{false}

\ifthenelse{\boolean{thesiscorrections}}{
  \newcommand{\corrections}[1]{\textcolor{red}{#1}}}
{
  \newcommand{\corrections}[1]{\textcolor{black}{#1}}
}

\ifthenelse{\boolean{changes}}{
  \newcommand{\changes}[1]{\textcolor{red}{#1}}}
{
  \newcommand{\changes}[1]{\textcolor{black}{#1}}
}

\ifthenelse{\boolean{examplechanges}}{
  \newcommand{\examplechanges}[1]{\textcolor{blue}{#1}}}
{
  \newcommand{\examplechanges}[1]{\textcolor{black}{#1}}
}

\ifthenelse{\boolean{refereecommentsvisible}}{
  \newcommand{\refereecommentsvisible}[1]{\textcolor{purple}{#1}}}
{
  \newcommand{\refereecommentsvisible}[1]{\textcolor{black}{#1}}
}

\begin{abstract}
    In this paper we present a framework for the construction and implementation of general virtual element spaces based on projections built from \refereecommentsvisible{constrained} least squares problems.
    Building on the triples used for finite element spaces, we introduce the concept of a VEM tuple which encodes the necessary building blocks to construct these projections.  
    Using this approach, a wide range of virtual element spaces can be defined.
    We discuss $H^k$-conforming spaces for $k=1,2$ as well as divergence and curl free spaces.
    This general framework has the advantage of being easily integrated into any existing finite element package and we demonstrate this within the open source software package \textsc{Dune}.
\end{abstract} 

    

\section{Introduction}\label{sec: intro}
In recent years research and development of the virtual element method (VEM) has skyrocketed. 
The method, which was introduced in \cite{beirao_da_veiga_basic_2013}, is an extension of the finite element method (FEM) and has many benefits. These include, but are not limited to, the handling of general polytopal meshes as well as the construction of spaces with additional structure such as arbitrary global regularity \cite{antonietti_conforming_2018,da2014virtual}. 
The method is highly flexible and as such has been applied to a wide range of problems.
The list is extensive and detailing all developments goes beyond the scope of this paper,
\changes{an overview can be found in \cite{Vem_and_Applications}.}

One of the main ideas behind the virtual element method is that the basis functions are considered virtual and do not need to be evaluated explicitly. 
In contrast to standard FEM, the local spaces may include functions which are not polynomials.   
Therefore, when setting up the stiffness matrix there are additional procedures which need to be carried out as projection operators need to be introduced and used when one of the entries of the stiffness matrix is a non-polynomial. 
It is important to stress that the computation of the projections is the only modification required to switch between a finite element and a virtual element discretisation.

A key part of the virtual element method is the construction of the aforementioned projection operators which are necessary to construct the discrete bilinear form. 
Usually, one projection operator is built which depends on the local contribution of the bilinear form \cite{beirao_da_veiga_basic_2013}.
This \changes{approach leads to some restrictions} as for example it does not extend easily to nonlinear problems or to PDEs with non-constant coefficients.
Also, it is not straightforward to integrate such a method into an existing finite element package which relies on the spaces not depending on the data of the PDE being solved.  
An alternative approach first introduced in \cite{ahmad_equivalent_2013} and extended to general second-order elliptic problems in \cite{cangiani_conforming_2015}, involves a VEM enhancement technique to ensure the computation of certain $L^2$-projections. 
This approach was extended to fourth-order problems in \cite{10.1093/imanum/drab003} and applied to a nonlinear example in \cite{dedner2021higher} where the starting point for the construction of projections is a \refereecommentsvisible{constrained} least squares (CLS) problem.

We aim to generalise this projection approach even further by introducing a general VEM framework
in which we describe how to define fully computable projections starting from what we call a generic \emph{VEM tuple}.  
The VEM tuple should be thought of as an extension of the well known \emph{finite element triple} \refereecommentsvisible{$(\element,\cP^\element,\Lambda^\element)$} 
\cite{clarlet1987finite}.
\changes{\refereecommentsvisible{Under the assumption of \emph{unisolvency} of $\Lambda^{\element} = \{ \lambda_1^{\element},\dots,\lambda_{N}^{\element} \}$,} this triple provides all the information required to construct the
nodal basis functions
\refereecommentsvisible{$\Phi^\element=\{ \phi_1^{\element}, \dots, \phi_N^{\element}\}$, i.e. a basis of $\cP^{\element}$ satisfying $\lambda_i^{\element}(\phi^{\element}_j)=\delta_{ij}$.}
Code to evaluate the basis
functions $\Phi^\element$ and their derivatives
for a given FEM triple forms the fundamental building block of most
finite element packages.}
In most cases changing the element type, the properties of the spaces, or even only the order of the method can lead to a considerable amount of new code. 
Implementing more specialised finite element spaces is often so cumbersome that few software packages provide these spaces. 
For example, spaces to handle
fourth-order problems in primal form are often not available except
sometimes in the form of the lowest order nonconforming Morley element
\cite{morley1968triangular}.

The introduction of a VEM tuple in this paper aims to extend the idea of a FEM triple as a fundamental building block for the implementation of a wide variety of VEM spaces. 
In particular, we incorporate all of the information required to build the projection operators in the VEM tuple.
One of the benefits of this generic approach is the ability to build spaces with additional properties.
We showcase this through examples detailing how to construct $H^k$-conforming for $k=1,2$, nonconforming, curl free and divergence free spaces suitable for e.g. acoustic vibration problems \cite{beirao2017virtual}, Stokes flow \cite{da_veiga_divergence_2015,dassi_bricks_2018}, to name but a few.

Extending the FEM concepts \corrections{i.e. replacing the evaluation of the nodal basis functions and their derivatives with the evaluation of projection operators,} reduces the complexity of adding virtual element methods to existing software frameworks.
\corrections{Such frameworks include for example} deal.II \cite{bangerth2007deal},
DUNE \cite{dunegridpaperII}, Feel++ \cite{prud2012feel}, FEniCS \cite{alnaes2015fenics}, Firedrake \cite{rathgeber2016firedrake}, or FreeFEM++ \cite{hecht2012new}.
This is in contrast to available VEM software implementations such as \cite{herrera2022numerical, balaje_kalyanaraman_2021_4561721, ortiz2019veamy, sutton2017virtual} which are often highly restrictive in what they offer.
Either only lowest order spaces are available or there is no functionality to solve nonlinear problems.
In this paper we explain the steps necessary to adapt a finite element implementation to include virtual element spaces. 
As proof of concept we demonstrate in two space dimensions how to integrate our general framework into a finite element package using the \textsc{Dune-Fem} module which is a part of the \textsc{Dune} software framework \cite{dedner2010generic}.
Switching between a finite element and a virtual element discretisation is  seamless for users due to the \textsc{Dune} Python frontend.  
To the best of our knowledge, this is the first implementation of the virtual element method within a large finite element software package to include spaces suitable for nonlinear problems in addition to spaces with the extra properties listed above.

The paper is organised as follows. 
In \Cref{sec: problem setup} we detail the notation and state the model problem followed by the details of the discretised problem.
The abstract virtual element framework is outlined in \Cref{sec: abstract framework} including the construction of the projection operators from a
given VEM tuple.
\changes{In \Cref{sec: examples} we provide a number of VEM tuples suitable for a wide range of different problems.
We restrict the main presentation to second-order problems in two space dimensions but provide an overview in \Cref{sec: fourth order problems} of how the framework extends to fourth-order problems and in \Cref{sec: ext to 3D} we sketch the extension of the approach to three space dimensions.
We outline how to use the general framework to add VEM spaces to existing FEM software packages in \Cref{sec: implementation details}}. 
Finally, numerical experiments are presented in \Cref{sec: numerics}. 
Note that for now our proof of concept implementation is restricted to two space dimensions.

\section{Problem setup}\label{sec: problem setup}
In this preliminary section we describe the notation and setup a model problem which we will use as an example during the discussion of the general concepts in the following.


\subsection{Mesh notation and assumptions}

Let $\mesh$ denote a tessellation of a domain $\Omega \subset \R^d$ for $d=2,3$, where $\Omega$ is a polygonal domain in 2D and a polyhedral domain in 3D. 
We assume the tessellation is made up of simple polygonal or polyhedral elements which do not overlap.
Moreover, we assume that the boundary of each element is made up of a uniformly bounded number of interfaces (edges in 2D and faces in 3D) where each boundary is either part of $\partial \Omega$ or shared with another element in $\mesh$.  
We denote an element in $\mesh$ by $\element$ with $h_{\element} := \text{diam}(\element)$, 
and let $h$ denote the maximum of the diameters over all elements in $\mesh$. 
We denote with $s$ a $(d-1)$-dimensional mesh interface, either an edge when $d=2$ or a face when $d=3$. The set of all interfaces in the mesh is denoted by $\cS_h$.
For $\element \in \mesh$ we denote by $\cS_h(\element)$ the set of all mesh interfaces which are part of $\element$, i.e., $\edgeSet := \{ s \in \cS_h : s \subset \partial \element \}$. 

We denote the space of polynomials over $\element \in \mesh$ with degree up to and including $\polOrder \in \N$ by $\prob_{\polOrder}(\element)$. 
We also require a hierarchical basis for $\prob_{\polOrder}(\element)$, which we denote by $\cM_{\polOrder}(\element)$. 
An example is given by scaled monomials which are defined as follows.
Firstly, for $k \in \N$, we denote by $\cM^{*}_{k}(\element)$ the scaled monomials of order $k$
\begin{align*}
    \cM^{*}_{k}(\element) = \left\{ \left( \frac{x-x_{\element}}{h_{\element}}\right)^{\alpha} \text{ with } \alpha \in \N^d, \, |\alpha|=k  \right\} 
\end{align*}
where $x_{\element}$ denotes the barycentre of $\element \in \mesh$, and $\alpha$ represents a multi-index with $|\alpha|=\sum_{i=1}^{d} \alpha_i$.
Then, we take as basis of $\prob_{\polOrder}(\element)$ the set $\cM_{\polOrder}(\element)$ where
\begin{align}\label{eqn: monomial basis}
    \cM_{\polOrder}(\element) := \cup_{k \leq \polOrder} \cM^{*}_{k}(\element).
\end{align}
Notice that we can also construct bases for polynomial spaces defined on interfaces $s$ in the same way denoted with $\cM_{\polOrder}(s)$. 
    

\subsection{Model problem and discrete setup}\label{sec: model problem}
\corrections{From now on, we restrict the presentation of the framework to two space dimensions, \changes{i.e. let} $\Omega \subset \R^2$ denote a polygonal domain.}
\corrections{We begin by setting up the model problem and describe the corresponding discrete problem.
\changes{Note that} we introduce the model problem with purposeful ambiguity to demonstrate the wide range of problems our framework can be applied to.
To this end, a}ssume that we have a form $a(\cdot,\cdot)$ defined on $V \times V$, where $V$ is the function space in which we define the associated PDE.
We assume $a(\cdot,\cdot)$ can be expressed in the following way  
\begin{align}\label{eqn: cts bilinear form}
    a(v,w) = \sum_{\element \in \mesh} a^\element(v,w)~,\; a^\element(v,w) = \int_{\element} D(v,\nabla v) \cdot \nabla w \, \dx + \int_{\element} m(v,\nabla v) w
    \, \dx  
\end{align}
for some possibly nonlinear $D,m$. 
\corrections{Further,} \refereecommentsvisible{for a linear functional $L$ defined on $V$,} we assume that there exists a unique solution $u \in V$ such that 
\begin{align}\label{eqn: cts problem}
    a(u,v) = \refereecommentsvisible{L(v)} \quad \forall v \in V.
\end{align}
\corrections{Since we introduce a model problem with deliberate ambiguity, alongside the abstract concepts, we also introduce a running example in the simplest setting.
That is, we make use of the two-dimensional Laplace problem in order to demonstrate the ideas behind our framework which are introduced throughout this section.}

\begin{example}[Conforming VEM example in 2D]\label{ex: 2D conforming VEM}
    \refereecommentsvisible{For a polygonal $\Omega \subset \R^2$,} consider the following two-dimensional Laplace problem \refereecommentsvisible{for some $f \in L^2(\Omega)$}:
    \begin{alignat*}{3}
        -\Delta u &= \refereecommentsvisible{f} , \quad &&\text{ in } \Omega,&&
    \end{alignat*}
    \changes{with \refereecommentsvisible{homogeneous} Dirichlet boundary conditions}.
    \corrections{Therefore, to fit into the above format,} we take the space ${V=H^1_0(\Omega)}$ and to obtain the bilinear form $a(\cdot,\cdot)$, we take \refereecommentsvisible{$D(u,\nabla u) = \nabla u$} and $m$ to be zero.
    This gives $a(u,v)=(\nabla u,\nabla v)$, \refereecommentsvisible{together with $L(v)=(f,v)$} it is clear that \eqref{eqn: cts problem} has a unique solution via the Lax-Milgram lemma.
\end{example}

\begin{remark}\label{rmk: vector valued}
    \changes{We point out that we can include vector valued problems with a form
    $a(\cdot,\cdot): \boldsymbol{V} \times \boldsymbol{V} \rightarrow \R$,
    where $\boldsymbol{V} := [V]^r$ for some $r>1$ and could for example be taken as
    $\boldsymbol{V} :=[H^1_0(\Omega)]^r$ in the case that $a(\boldsymbol{v},\boldsymbol{w}) :=(\nabla \boldsymbol{v}, \nabla \boldsymbol{w})$.}
\end{remark}


\corrections{Next,} given a space $\widetilde{V}$ not necessarily satisfying $\widetilde{V} \subset V$, we define a discrete bilinear form $a_h : \widetilde{V} \times \widetilde{V} \rightarrow \R$ such that for any $v_h, w_h \in \widetilde{V}$
\begin{align*}
    a_h(v_h,w_h) = \sum_{\element \in \mesh} a_h^\element (v_h,w_h)
\end{align*}
where $a_h^\element : \widetilde{V}|_{\element} \times \widetilde{V}|_{\element} \rightarrow \R$ is the restriction of the bilinear form $a_h$ to an element $\element$. 
Denote this restriction of $\widetilde{V}$ to an element $\element$ by $\widetilde{V}^\element$. 
We assume further that the discrete form has the following form
\begin{equation}\label{eqn: discrete bilinear form}
    \begin{split}
        a_h^\element (v_h,w_h) = 
        &\int_{\element} D(\valueProj v_h, \gradProj v_h ) \cdot \gradProj w_h \, \dx 
        + \int_{\element} m(\valueProj v_h, \gradProj v_h) \valueProj w_h \, \dx
        \\
        &+ \corrections{( \bar D + \bar m h_{\element}^{2} )} S^\element(v_h-\valueProj v_h, w_h - \valueProj w_h)
    \end{split}
\end{equation}
for any $v_h, w_h \in \widetilde{V}^{\element}$, where $D,m$ have been given in \cref{eqn: cts bilinear form}, \corrections{and $\bar D$, $\bar m$ are constants necessary to achieve the right scaling of the so called VEM \emph{stabilisation} form $S^{\element}(\cdot,\cdot)$. 
$S^\element(\cdot,\cdot)$ denotes a symmetric, positive definite bilinear form necessary to ensure the coercivity of $a_h^\element(\cdot,\cdot)$.}
\corrections{Furthermore,} $\Pi^{\element}_{\mu}$ denotes a value projection for $\mu=0$ and a gradient projection for $\mu=1$ chosen so that 
\begin{align}
    \Pi^{\element}_{\mu} v_h \sim \nabla^{\mu} v_h |_{\element}, \quad \mu = 0,1.
    \label{eqn: projections}
\end{align}

\refereecommentsvisible{In order to discretise the right hand side of \eqref{eqn: cts problem}, for any $v_h \in \widetilde{V}$ we define $L_h(v_h) := \sum_{E \in \mesh} L_h^{\element}(v_h)$.
In the case of \Cref{ex: 2D conforming VEM}, we define $L^{\element}_h (v_h): = ( f, \valueProj v_h)_{\element}$, i.e. we utilise the value projection. 
Note that for the remainder of this paper we focus on the bilinear form.
In particular, the definitions of the projections will be made precise in the virtual element setting in \Cref{sec: abstract framework}.}
\changes{
\begin{remark}
For a classical conforming FEM discretisation we would simply take $S^{\element}=0$ and
$\Pi^{\element}_{\mu} v_h = \nabla^{\mu} v_h |_{\element}$ so that
$a_h^\element (v_h,w_h) = \refereecommentsvisible{a^{\element}} (v_h,w_h)$ \refereecommentsvisible{and $L_h^{\element}(v_h)=L^{\element}(v_h)=(f,v_h)_{\element}$.} 
\end{remark}
}  
\begin{definition}[Degrees of freedom]\label{def: degrees of freedom}
    For $\element \in \mesh$, we define a set of \emph{degrees of freedom} \refereecommentsvisible{(dof set)} $\elementDofSet$ as a set of functionals $\lambda : \widetilde{V} \rightarrow \R$. 

    For a given dof set $\elementDofSet$ we define the subset $\Lambda^{\element,s}$
    for $s \in \cS_h(\element)$ as follows
    \begin{align}\label{eqn: interface dof set}
        \Lambda^{\element,s} &:= \{ \lambda \in \elementDofSet : \ \lambda(v^{\element} ) = \lambda(w^{\element}) \text{ if } v^{\element} |_{\bar s} = w^{\element} |_{\bar s} \}.
    \end{align}
\end{definition}

\begin{assumption}\label{assumption: dofs}
    We assume the following. 
    \begin{enumerate}[label=(\textit{A}\arabic*)]
        \item The dof set $\elementDofSet$ is unisolvent, i.e., a function $v_h$ in the local discrete space is uniquely determined by its degrees of freedom \refereecommentsvisible{(dofs)}.
        \item The dofs depend only on values of $v_h$ but not on derivatives. 
        \label{assumption: A2}
        \item The polynomial space $\prob_{\polOrder}(\element)$ is a subset of the local discrete VEM space, \refereecommentsvisible{where $\polOrder$ represents a positive integer}. 
    \end{enumerate}
\end{assumption}
Note that \ref{assumption: A2} will be generalised in a later section (\Cref{sec: fourth order problems}) \corrections{when we consider fourth-order problems}.
We \corrections{now} give an example set of dofs for \Cref{ex: 2D conforming VEM}. 

\begin{continuance}{ex: 2D conforming VEM}\label{2.1 dofs}
    As considered in e.g. \cite{beirao_da_veiga_basic_2013,cangiani_conforming_2015}, for a polygon $\element \in \mesh$, we consider the set of dofs $\Lambda^{\element}$ for the $H^1$-conforming VEM.
    For $\polOrder \geq 1$, we define the following dofs.
    \begin{enumerate}[label=(\alph*)]
        \item 
        \begin{itemize}
            \item The values of $v_h$ for each vertex of $\element$.
            \item For $\polOrder>1$, the moments of $v_h$ up to order $\polOrder-2$ on each edge $s \in \edgeSet$,
            \begin{align}\label{eqn: h1 conforming edge dofs}
                \frac{1}{|s|} \int_s v_h m_{\polOrder-2} \, \ds \quad &\forall m_{\polOrder-2} \in \cM_{\polOrder-2}(s).
            \end{align} 
        \end{itemize}
        \label{item: vertex and edges dofs c0 example}
        \item For $\polOrder > 1$, the moments of $v_h$ up to order $\polOrder-2$ inside the element $\element$
        \begin{align*}
            \frac{1}{|\element|} \int_{\element} v_h m_{\polOrder-2} \, \dx \quad &\forall m_{\polOrder-2} \in \cM_{\polOrder-2}(\element). 
        \end{align*}
        \label{eqn: h1 conforming inner dofs}
    \end{enumerate}
    In this example, we have $N^s + N^s(\polOrder-1) + \half \polOrder(\polOrder-1)$ total dofs where $N^s$ denotes the number of edges in the polygon $\element$. 
    The dof set $\Lambda^{\element,s}$ contains the \corrections{dofs in \ref{item: vertex and edges dofs c0 example}, that is,} the two dofs $v_h(s^{\pm})$ where $s^{\pm}$ denotes the vertices attached to the edge $s$, as well as the $|\cM_{\polOrder-2}(s) | = \polOrder-1$ edge dofs.
\end{continuance}

\corrections{Key to the virtual element discretisation is the concept of computability which we capture in the next definition.}

\begin{definition}[Computable]
    Given $F \in \mesh$ (or ${F \in \edgeSet}$ for some $\element \in \mesh$), we say that a quantity is \emph{computable from $\Lambda^{F}$} (or $\Lambda^{\element,s}$) if it is a linear combination of $\lambda \in \Lambda^{F} (\text{ or } \Lambda^{\element,s})$.
\end{definition}

\section{Abstract virtual element framework}\label{sec: abstract framework}
In this section we introduce an approach to construct general VEM spaces based on ideas borrowed from the abstract FEM framework where a triple is used to define the local space.
We begin by defining the VEM tuple before describing the projection framework \refereecommentsvisible{in two space dimensions}.
\refereecommentsvisible{Recall that the projections are necessary to construct the discrete form \eqref{eqn: discrete bilinear form}. In particular, the element value projection $\valueProj$ and gradient projection $\gradProj$ are chosen so that they are as ``close'' as possible to the true value and gradient operators \eqref{eqn: projections}.}

Within the projection framework we introduce two value projections, one for each element $\element \in \mesh$ and one for \corrections{each edge $s \in \edgeSet$.}
These \refereecommentsvisible{value} projections are defined as solutions of a \refereecommentsvisible{constrained} least squares problem with the constraint set specified by the VEM tuple.
\changes{The projections for the higher order derivatives e.g. the gradient or hessian projections are then defined hierarchically using these projections, independent of the space.
We discuss the gradient projection here since we are focusing on second-order problems and
delay the description of the generic construction of the hessian projection to \Cref{sec: fourth order problems}.}


\subsection{Virtual element tuple}\label{sec: VEM tuple}
The starting point for the general framework is the notion of a virtual element tuple.
The tuple contains all of the building blocks for the computation of the projection operators.

\begin{definition}[Virtual element tuple]\label{def: vem tuple}
    We define a \emph{virtual element tuple} $\cV^{F}$
    \begin{align}\label{eqn: VEM tuple}
        \cV^{F} = {(F, \valueBasisSet^F, \Lambda^F, \cC^{F}_0, \gradBasisSet^F)}
    \end{align}
    consisting of the following.
    \begin{itemize}
        \item A mesh object $F$ with either $F \in \mesh$ or $F \in \cS_h$.
        \item A basis set $\valueBasisSet^F$ used for a value projection e.g. a subset of \changes{$[\cM_{\polOrder}(F)]^r$} (see \cref{eqn: monomial basis}),
        \examplechanges{where $r=1$ for scalar and $r>1$ for vector valued spaces, $|\valueBasisSet^F| \leq |\Lambda^F|$, and  when $F \in \mesh$, $\valueBasisSet^{F}$ is a subset of the local discrete VEM space.}
        \item A set of dofs $\Lambda^F$ subject to the conditions in \cref{assumption: dofs}.
        \item A set of linear functionals $\cC^{F}_{0}$ computable from $\Lambda^{F}$. These will be called constraints in the following.
        \item When $F \in \mesh$, a basis set $\gradBasisSet^F$ used for a gradient projection. 
        \corrections{Note that since the basis set $\gradBasisSet^F$ is used to construct the \emph{gradient} projection, it must at least be a} \changes{vector valued} quantity, \changes{e.g. a subset of $[\cM_{\polOrder-1}(F)]^{r \times 2}$}.
    \end{itemize}
    For the set of VEM tuples $(\cV^{\element})_{\element \in \mesh}$ and $(\cV^{s})_{s \in \cS_h}$ we say that the VEM tuples are \emph{compatible} if $\Lambda^{s} \subset \Lambda^{\element}$ whenever $s \in \edgeSet$.
\end{definition}

\begin{continuance}{ex: 2D conforming VEM}\label{basis sets}
    \examplechanges{For our conforming VEM running example,
    the dof set $\elementDofSet$ is taken as those described previously in \Cref*{ex: 2D conforming VEM}\ref{item: vertex and edges dofs c0 example}-\ref{eqn: h1 conforming inner dofs}.
    For each $\element \in \mesh$, we take as basis set the monomial basis $\valueBasisSet^{\element}$ from \cref{eqn: monomial basis} of order $\polOrder$ and we take the basis for the gradient projection, $\gradBasisSet^{\element}$, to be the monomials of order $\polOrder-1$, i.e. $\gradBasisSet^{\element} := [\cM_{\polOrder-1}(\element)]^2$.
    Also, for this example, for each $s \in \edgeSet$, we take $\valueBasisSet^s$ to be $\cM_{\polOrder}(s)$. That is, we shall define the value projections such that $\valueProj v_h \in \prob_{\polOrder}(\element)$ and $\faceProj v_h \in \prob_{\polOrder}(s)$. Furthermore, we shall define the gradient projection such that $\gradProj v_h \in [\prob_{\polOrder-1}(\element)]^2$.
    The dof set for the edge tuple is taken to be $\Lambda^s := \Lambda^{\element,s}$.
    We will discuss the constraints sets $\cC^{\element}_0,\cC^s_0$ for this example later in this section.
    }
\end{continuance}

For the remainder of this section we now fix $\element \in \mesh$ and $\Lambda^{\element}$.
In order to ensure compatibility of the VEM tuples we take $\Lambda^s := \Lambda^{\element,s}$ which is defined in \cref{eqn: interface dof set} for $s \in \edgeSet$. Then it is clear that $\Lambda^s \subset \Lambda^{\element}$.
\begin{remark}\label{remark: vem tuples}
    For second-order elliptic problems, when $F \in \edgeSet$ the basis set $\gradBasisSet^F$ is not prescribed as we do not require higher order derivative projections on the edges, this is due to \ref{assumption: A2}.
    The more general case will be discussed later.
\end{remark}

\begin{remark}
\examplechanges{
    To construct finite element spaces a triple of the form $(F, \valueBasisSet^F, \Lambda^F)$ is commonly used. The main difference to the VEM tuple is that in the FEM setting $\valueBasisSet^F$ is a basis of the local space, in fact usually the local space itself is used in the FEM triple.
    In our VEM framework, we can have $|\valueBasisSet^{\element}| < N^{\element}$ where $N^{\element}$ denotes the dimension of the discrete space. So only a (mostly polynomial) subspace of the local space is specified by the VEM tuple.
    }
\end{remark}

Both the basis sets $\cB^F_{\mu}$ for $\mu=0,1$ and the constraint set $\cC^F_0$  are kept ambiguous to allow for the construction of VEM spaces with different properties.
In Section~\ref{sec: examples} we give concrete examples for these choices, but in general the basis sets will be sets of monomials used as a basis for various polynomial spaces.


\subsection{Gradient projection}
We now give the details of the projection framework. 
Given an element $\element \in \mesh$, attached \corrections{edges}, and VEM tuples $\cV^{\element}$, $\cV^{s}$ (see \Cref{def: vem tuple}) we introduce value projections $\Pi^{\element}_0,$ $\edgeProj$ for each \corrections{edge} $s$, and a gradient projection $\gradProj$. 
In order to illustrate the hierarchical construction process, we start by giving the definition of the gradient projection \corrections{in terms of the value projections $\valueProj$, $\edgeProj$ to be determined later.}
The definition of the gradient projection comes from an application of integration by parts followed by a replacement of the lower order terms with computable projections, as shown in the next definition.
This ensures that the gradient projection is fully computable.
\changes{Note that it is more usual in the VEM literature to first define the gradient projection and then to use that to define the value projection.}

\begin{definition}[Gradient projection]\label{def: gradient projection}
    For a discrete function $v_h$ and an element $\element \in \mesh$,
    we define the \emph{gradient projection} $\gradProj$ into the span of the basis set $\gradBasisSet^{\element}$, 
    where $\gradBasisSet^{\element}$ is provided by the VEM tuple $\cV^{\element}$ \eqref{eqn: VEM tuple}.
    Then \corrections{for any $v_h \in \widetilde V^{\element}$,} $\gradProj$ \changes{satisfies} 
    \begin{align*}\label{eqn: grad proj}
        \int_{\element} \gradProj v_h \cdot q \, \dx = - \int_{\element} \valueProj v_h \nabla \cdot q \, \dx + \sum_{s \in \cS_h(\element)} \int_s \faceProj v_h (n \cdot q) \, \ds
    \end{align*}
    for all $q \in \gradBasisSet^{\element}$, 
    where $n$ denotes the outward pointing normal to the \corrections{edge} $s$.
\end{definition}
Motivated by this definition, it is clear why we need to construct an inner value projection and projections on each \corrections{edge} $s \in \cS_h(\element)$.

\subsection{Value projections}
In this section, we now let $F$ denote either the element $\element$ or an \corrections{edge} $s \in \cS_h(\element)$.
We stress that there is no difference in the construction of the value projection whether $F$ is an element or an \corrections{edge}. 
Either way, these projections are defined using a \refereecommentsvisible{constrained} least squares problem, which we assume is uniquely solvable.
\changes{Details on how to implement this projection can be found in \Cref{sec: implementation details} where we also discuss the unique solvability of the \refereecommentsvisible{constrained} least squares problem.}


\begin{definition}[Value projection]\label{def: value projection}
    For a discrete function $v_h$, we define the \emph{value projection} $\meshProj$ into the span of the basis set $\valueBasisSet^F$, 
    where $\valueBasisSet^F$ is provided by $\cV^F$.
    We define the value projection $\meshProj$ as the solution to the following CLS problem 
    \begin{align*}
        {\rm min} \sum_{\lambda \in \Lambda^{F}} (\lambda (\meshProj v_h - v_h))^2, \qquad \text{subject to} \quad \cC ( \meshProj v_h - v_h ) = 0 \quad \forall \cC \in \cC^F_0.
    \end{align*}
\end{definition}

We revisit \cref{ex: 2D conforming VEM} again, giving concrete examples for the constraint set choices in the VEM tuples for this example. 

\begin{continuance}{ex: 2D conforming VEM}\label{constraint sets}
    \examplechanges{We have so far defined all components of the VEM tuple except the constraint sets. One important aim of the constraints is to make sure that the value projection is an $L^2$ projection into as large a polynomial set as possible. Consequently, if for some polynomial $q$ the right hand side $\int_\element v_hq$ of the $L^2$ projection of a discrete function $v_h$ is computable from its degrees of freedom, then the value projection should be constrained to satisfy $\int_\element\meshProj v_hq = \int_\element v_hq$. In the case of the $H^1$-conforming space with dofs given by \Cref*{ex: 2D conforming VEM}\ref{item: vertex and edges dofs c0 example}-\ref{eqn: h1 conforming inner dofs}, we therefore take as constraint set a scalar multiple of the inner moments:
    \begin{align}\label{eqn: h1 conf value constraints}
        \cC (v_h) := \int_{\element} v_h m_{\polOrder-2} \, \dx, \quad m_{\polOrder-2} \in \cM_{\polOrder-2}(\element).
    \end{align}
    The constraint set for the edges is given by the edge dofs only from \cref{eqn: h1 conforming edge dofs}.  
    This leads to the following CLS problem for the value projection on each edge with a unique solution in $\prob_{\polOrder}(s)$.
    \begin{align}
        {\rm min} \sum_{\lambda \in \interfaceDofSet} (\lambda (\faceProj v_h - v_h))^2 \quad \text{subject to} 
        \int_s \left( \faceProj v_h - v_h \right) m_{\polOrder-2} \, \ds = 0, \ m_{\polOrder-2} \in \cM_{\polOrder-2}(s).
        \label{eqn: h1 conforming edge constraints}
    \end{align}
    }
\end{continuance}

\begin{remark}
    \refereecommentsvisible{
    We note that for the upcoming $H^1$-conforming and nonconforming examples we can show using the standard approach of ``extended VEM spaces'' as in \cite{cangiani_conforming_2015,10.1093/imanum/drab003} that for any discrete $v_h$ and for $\mu=0,1$
    \begin{align*}
        \Pi^{\element}_{\mu} v_h = \cP^{\element}_{\polOrder-\mu} (\nabla^{\mu} v_h).
    \end{align*}
    This says that the value projection is indeed the $L^2$-projection, denoted by $\cP^{\element}_{p}$ for $p\in\N$, and the gradient projection is the $L^2$-projection of the gradient.
    This follows due to our choice of constraints, combined with the standard ``extended VEM space'' approach.
    We therefore obtain the usual $L^2$ property necessary for convergence analysis.
    }
\end{remark}
\subsection{Examples}\label{sec: examples}
In this section we give concrete examples of VEM tuples for specific problems. 
Our generic framework allows us to build further VEM spaces with additional properties and a unified way to construct the projection operators for these problems.
It is important to reiterate that assuming we are given $\element \in \mesh$ and  $\elementDofSet$, when $s \in \edgeSet$ we take $\Lambda^s := \Lambda^{\element,s}$ in the VEM tuple $\cV^s$.
Recall that $\Lambda^{\element,s}$ is defined in \cref{eqn: interface dof set} and satisfies $\Lambda^{\element,s} \subset \elementDofSet$ thus ensuring compatibility of the VEM tuples.

\subsubsection{Example: \texorpdfstring{$H^1$}{H1}-conforming VEM}\label{sec: h1 VEM 2D}
We \changes{summarise} the VEM tuple for the $H^1$-conforming VEM in two dimensions from the running Example~\ref{constraint sets} in \cref{tab: h1 conf 2d}.
\begin{table}[htbp]
    \footnotesize
    \caption{VEM tuple summary for the $H^1$-conforming VEM in two dimensions. \changes{We introduce $q \in \N$ with $q \leq \polOrder$ for the gradient basis set $\gradBasisSet^{\element}$.}}\label{tab: h1 conf 2d}
    \begin{center}
      \begin{tabular}{|c|c|c|c||c|c|} \hline
       $\valueBasisSet^{\element}$ & $\elementDofSet$ & $\cC^{\element}_0$ & $\gradBasisSet^{\element}$ & $\valueBasisSet^{s}$ & $\cC^s_0$ \\ 
       \hline
        $\cM_{\polOrder}(\element)$ &
        \Cref*{ex: 2D conforming VEM}\ref{item: vertex and edges dofs c0 example}-\ref{eqn: h1 conforming inner dofs} &
        \cref{eqn: h1 conf value constraints} &
        $[\cM_{q}(\element)]^2$ &
        $\cM_{\polOrder}(s)$ &
        \cref{eqn: h1 conforming edge constraints} 
        \\
        \hline
      \end{tabular}
    \end{center}
\end{table}

\begin{remark}[Stabilisation free spaces]\label{rmk: stab free}
    \changes{As discussed in Example~\ref{basis sets}, the natural choice for $q$ in \Cref{tab: h1 conf 2d} for the gradient basis set $\gradBasisSet^{\element}$ is $q:=\polOrder-1$ however, we can take $q:=\polOrder$ following ideas for stabilisation free VEM as introduced for the lowest order space in \cite{berrone2021lowest}. }
\end{remark}


\begin{remark}\label{rmk: smaller basis sets}
    Another choice of basis set for the gradient projection which is smaller than $[\cM_{\polOrder-1}(\element)]^2$ could be $\gradBasisSet^{\element} := \nabla \cM_{\polOrder}(\element)$.
    This choice of basis set is suitable for constant coefficient problems as otherwise it leads to suboptimal convergence rates for varying coefficients as remarked in \cite{beirao2016virtual}.
\end{remark}
\begin{remark}[Serendipity spaces]\label{rem: serendipity VEM}
\examplechanges{
The serendipity approach discussed in \cite{da_veiga_serendipity_2015} can be easily incorporated by choosing fewer moment degrees of freedom in the interior of each element without needing a change to the projection operators.
}
\end{remark}


\subsubsection{Example: \texorpdfstring{$H^1$}{H1}-nonconforming VEM}\label{sec: ex h1 nonconforming VEM}
In this example we consider a nonconforming VEM space suitable for second-order problems taking the degrees of freedom considered in e.g. \cite{cangiani_conforming_2015, chen2020nonconforming,de_dios_nonconforming_2014}. 
\begin{definition}\label{def: h1 nonconforming dofs}
    For $\polOrder \geq 1$, we define the following dofs.
    \begin{itemize}
        \item The moments of $v_h$ up to order $\polOrder-1$ on each \corrections{edge} $s \in \edgeSet$
        \begin{align}\label{eqn: h1 nonconforming edge dofs}
            \frac{1}{|s|} \int_s v_h m_{\polOrder-1} \, \ds \quad &\forall m_{\polOrder-1} \in \cM_{\polOrder-1}(s).
        \end{align}
        \item For $\polOrder >1$, the moments of $v_h$ up to order $\polOrder-2$ inside the element $\element$
        \begin{align}\label{eqn: h1 nonconforming inner dofs}
            \frac{1}{|\element|} \int_{\element} v_h m_{\polOrder-2} \, \dx \quad &\forall m_{\polOrder-2} \in \cM_{\polOrder-2}(\element).
        \end{align}
    \end{itemize}
\end{definition}
Note that the set $\Lambda^{\element,s} \subset \elementDofSet$ includes the dofs in \cref{eqn: h1 nonconforming edge dofs} for each $s \in \cS_h(\element)$. 
The constraints for the edge value projection are a scalar multiple of the edge moments described in \cref{eqn: h1 nonconforming edge dofs}. Each constraint is therefore of the form 
\begin{align}\label{eqn: h1 nonconf edge constraints}
    \cC(v_h):= \int_{s} v_h m_{\polOrder-1} \, \ds, \quad m_{\polOrder-1} \in \cM_{\polOrder-1}(s).
\end{align}
\Cref{tab: h1 nonconforming vem tuples} summarises the choices for the other quantities of the VEM tuple.
\begin{table}[htbp]
    \footnotesize
    \caption{VEM tuples for the nonconforming VEM. In particular, we define the value projection into $\text{Span}(\cM_{\polOrder}(\element))=\prob_{\polOrder}(\element)$ and the \corrections{edge} projections into $\text{Span}(\cM_{\polOrder-1}(s))=\prob_{\polOrder-1}(s)$.}\label{tab: h1 nonconforming vem tuples}
    \begin{center}
      \begin{tabular}{|c|c|c|c||c|c|} \hline 
       $\valueBasisSet^{\element}$ & $\elementDofSet$ & $\cC^{\element}_0$ & $\gradBasisSet^{\element}$ & $\valueBasisSet^{s}$ & $\cC^s_0$ \\ 
       \hline
        $\cM_{\polOrder}(\element)$ & 
        \cref{def: h1 nonconforming dofs} & 
        \cref{eqn: h1 conf value constraints} & 
        $[\cM_{q}(\element)]^2$ & $\cM_{\polOrder-1}(s)$ & 
        \cref{eqn: h1 nonconf edge constraints} \\
        \hline
      \end{tabular}
    \end{center}
\end{table}
\changes{Similarly to \cref{rmk: stab free}, the natural choice is $q:=\polOrder-1$ however in the spirit of \cite{berrone2021lowest} we can also take $q:=\polOrder$ for a stabilisation free VEM for this nonconforming example.}

\begin{remark}
    \changes{The main difference to the conforming case is that $\faceProj$ is only into $\prob_{\polOrder-1}(s)$.
    This is sufficient for computing the gradient projection which only requires moments
    in $\gradBasisSet^{\element} = [\cM_{\polOrder-1}(\element)]^2$ to be exact}.
\end{remark}


\subsubsection{Example: Divergence free VEM}\label{sec: div free}
For this example we present the projection approach in the divergence free setting based on the \emph{reduced} dof set described in \cite{da_veiga_divergence_2015}.
\begin{definition}\label{def: div free reduced dofs}
    For $\polOrder \geq 2$, we define the following dofs for a discrete vector field
    $\boldsymbol{v}_h$.
    \begin{itemize}
        \item The values of $\boldsymbol{v}_h$ for each vertex of $\element$.
        \item The moments up to order $\polOrder-2$ for each \corrections{edge} $s \in \edgeSet$,
        \begin{align*}
            \frac{1}{|s|} \int_s \boldsymbol{v}_h \cdot \boldsymbol{m}_{\polOrder-2} \, \ds \quad &\forall \boldsymbol{m}_{\polOrder-2} \in [\cM_{\polOrder-2}(s)]^2.
        \end{align*}
        \item The moments of $\boldsymbol{v}_h$ inside the element $\element$
        \begin{align*}
            \frac{1}{|\element|} \int_{\element} \boldsymbol{v}_h \cdot \boldsymbol{m}_{\polOrder-2}^{\perp} \, \dx \quad &\forall \boldsymbol{m}_{\polOrder-2}^{\perp} \in \boldsymbol{x}^{\perp}[\cM_{\polOrder-3}(\element)],
        \end{align*}
    \end{itemize}
    with the notation $\boldsymbol{x}^{\perp} := (x_2,-x_1)$. 
\end{definition}

\corrections{For this example we also assume as in \cite{da_veiga_divergence_2015} that each discrete function in the local VEM space satisfies $\nabla \cdot \boldsymbol{v}_h \in \prob_{0}(\element)$ and $\boldsymbol{v}_h|_{s} \in [\prob_{\polOrder}(s)]^2$ for each $s \in \edgeSet$.
This allows us to show important properties of the projections.}

Up until now we have always taken the element constraint set $\cC_0^{\element}$ to contain the inner dofs. However, in this example we take some additional constraints as well as the inner dofs.
\begin{align}
    \cC(\boldsymbol{v}_h) &:= \int_{\element} \boldsymbol{v}_h \cdot \boldsymbol{m}_{\polOrder-2}^{\perp} \, \dx, \quad \boldsymbol{m}_{\polOrder-2}^{\perp} \in \boldsymbol{x}^{\perp}[\cM_{\polOrder-3}(\element)],
    \label{eqn: div free element constraints}
\intertext{and in addition}
    \cC(\boldsymbol{v}_h) &:= \int_{\element} \boldsymbol{v}_h \cdot \nabla m_{\polOrder-1} \, \dx, \quad m_{\polOrder-1} \in \cM_{\polOrder-1}(\element) \backslash \cM_{0}(\element).
    \label{eqn: div free extra constraints}
\end{align}
\corrections{The constraints for the edge value projection are the vector version of those used in the $H^1$-conforming example and so we are able to define $\edgeProj \boldsymbol{v}_h \in [\prob_{\polOrder}(s)]^2$. Since we assume $\boldsymbol{v}_h|_{s} \in [\prob_{\polOrder}(s)]^2$, it follows that $\edgeProj \boldsymbol{v}_h = \boldsymbol{v}_h |_{s}$.}

Note that the constraints in \cref{eqn: div free extra constraints} are indeed computable \corrections{and necessary to show properties of the discrete space - see \cref{rmk: trace of grad proj}.}
\corrections{Computability} follows from applying integration by parts to the right hand side, 
\begin{align*}
    \int_{\element} \boldsymbol{v}_h \cdot \nabla m_{\polOrder-1} \, \dx 
    &= 
    - \int_{\element} \nabla \cdot \boldsymbol{v}_h m_{\polOrder-1} \, \dx + \sum_{s \in \cS_h(\element)} \int_{s} (\boldsymbol{v}_h \cdot \boldsymbol{n}) m_{\polOrder-1} \, \ds.
\end{align*}
Assuming that $m_{\polOrder-1}$ is in an orthonormal basis, the first term disappears due to the assumption $\nabla \cdot \boldsymbol{v}_h \in \prob_{0}(\element)$.
The boundary term can also be computed using the edge projection since $\boldsymbol{v}_h|_s \in [\prob_{\polOrder}(s)]^2$. 

We now have all the ingredients to provide the VEM tuple given in \cref{tab: div free vem}.
\begin{table}[htbp]
    \footnotesize
    \caption{VEM tuples for the divergence free VEM. Notice that this example is vector valued and as such the basis set for the element value projection is $[\cM_{\polOrder}(\element)]^2$. }\label{tab: div free vem}
    \begin{center}
      \begin{tabular}{|c|c|c|c||c|c|} \hline 
       $\valueBasisSet^{\element}$ & $\elementDofSet$ & $\cC^{\element}_0$ & $\gradBasisSet^{\element}$ & $\valueBasisSet^{s}$ & $\cC^s_0$ \\ 
       \hline
        $[\cM_{\polOrder}(\element)]^2$ & 
        \cref{def: div free reduced dofs} & 
        \cref{eqn: div free element constraints}-\cref{eqn: div free extra constraints} & 
        $[\cM_{\polOrder-1}(\element)]^{2 \times 2}$ &
        $[\cM_{\polOrder}(s)]^2$ & 
        \cref{eqn: h1 conforming edge constraints} \\
        \hline
      \end{tabular}
    \end{center}
\end{table}

\begin{remark}\label{rmk: trace of grad proj}
    Using integration by parts, the constraints in \cref{eqn: div free extra constraints}, and the exactness of the edge projection, we can show that $\text{tr} (\gradProj v_h) \in \prob_{0}(\element)$.
    In particular, for any test function $m_{\beta} \in \cM_{\polOrder-1}(\element) \backslash \cM_{0}(\element)$ it holds that
    \begin{align*}
        \int_{\element} \text{tr} (\gradProj \boldsymbol{v}_h) m_{\beta} \, \dx
        &= - \int_{\element} \valueProj \boldsymbol{v}_h \cdot \nabla m_{\beta} \, \dx + \sum_{s \in \edgeSet} \int_s ( \faceProj \boldsymbol{v}_h \cdot \boldsymbol{n}) m_{\beta}
        \, \ds 
        \\
        &= \int_{\element} (\nabla \cdot \boldsymbol{v}_h ) m_{\beta} \, \dx.
    \end{align*}
    From an implementation perspective, this property is highly valuable as it \refereecommentsvisible{is one possibile approach which} avoids the construction of additional projection operators for the divergence. 
    We implement the divergence directly as the trace of the gradient projection
    which is likely to be the default implementation available in a finite-element code.
\end{remark}


\subsubsection{Example: Curl free VEM}\label{sec: curl free}
We now show how our generic approach can be applied to construct a curl free space using the dof set considered for an acoustic vibration problem in \cite{beirao2017virtual}. 
\begin{definition}\label{def: curl free dofs}
    For $\polOrder \geq 0$, we define dofs for a vector field $\boldsymbol{v}_h$:
    \begin{itemize}
        \item The moments up to order $\polOrder$ for each \corrections{edge} $s \in \edgeSet$,
        \begin{align}\label{eqn: curl free edge dofs}
            \frac{1}{|s|} \int_s (\boldsymbol{v}_h \cdot \boldsymbol{n}) m_{\polOrder} \, \ds \quad &\forall m_{\polOrder} \in \cM_{\polOrder}(s).
        \end{align}
        \item The moments of $\boldsymbol{v}_h$ inside the element $\element$,
        \begin{align}
            \frac{1}{\sqrt{|\element|}} \int_{\element} \boldsymbol{v}_h \cdot \nabla m_{\polOrder} \, \dx \quad &\forall m_{\polOrder} \in \cM_{\polOrder}(\element) \backslash \cM_0(\element). 
        \end{align}
    \end{itemize}
\end{definition}

\corrections{As in \cite{beirao2017virtual}, for this example we assume further that discrete functions in the local VEM space satisfy $\nabla \cdot \boldsymbol{v}_h \in \prob_{\polOrder}(\element)$ and $(\boldsymbol{v}_h \cdot \boldsymbol{n}) \in \prob_{\polOrder}(s)$ for each $s \in \edgeSet$.}

As in the previous example, we take as constraints the inner dofs as well as additional constraints based on the additional property $\nabla \cdot v_h \in \prob_{\polOrder}(\element)$. They are of the following form.
\begin{align}
    \cC(\boldsymbol{v}_h) &:= \int_{\element} \boldsymbol{v}_h \cdot \nabla m_{\polOrder} \, \dx, &&m_{\polOrder} \in \cM_{\polOrder}(\element),
    \label{eqn: curl free constraints}
    \\
    \cC(\boldsymbol{v}_h) &:= \int_{\element} \boldsymbol{v}_h \cdot \nabla m_{\polOrder+1} \, \dx, &&m_{\polOrder+1} \in \cM_{\polOrder+1}(\element) \backslash \cM_{\polOrder}(\element).
    \label{eqn: extra curl free constraint}
\end{align}
The constraints for the edge value projection are again a scalar multiple of the edge dofs in \cref{eqn: curl free edge dofs} and, for any $s \in \cS_h(\element)$, are of the form 
\begin{align}\label{eqn: edge curl free constraints}
    \changes{\int_s \boldsymbol{v}_h \cdot \boldsymbol{m_{\polOrder}} \, \ds,
        \quad \boldsymbol{m_{\polOrder}} \in \boldsymbol{n}\cM_{\polOrder}(s).}
\end{align}
These constraints allows us to define the edge projection into $\prob_{\polOrder}(s)\boldsymbol{n}$ and therefore $\edgeProj \boldsymbol{v}_h = \boldsymbol{v}_h |_{s}$.

It is clear that the constraints in \cref{eqn: curl free constraints} are computable using the dofs whereas the constraints in \cref{eqn: extra curl free constraint} reduce to the following 
\begin{align*}
    \int_{\element} \boldsymbol{v}_h \cdot \nabla m_{\polOrder+1} \, \dx
    &= 
    - \int_{\element} \nabla \cdot \boldsymbol{v}_h m_{\polOrder+1} \, \dx + \sum_{s \in \cS_h(\element)} \int_s (\boldsymbol{v}_h \cdot \boldsymbol{n}) m_{\polOrder+1} \, \ds.
\end{align*}
Assuming that we have an orthonormal basis for $\prob_{\polOrder+1}(\element)$, the first term vanishes since we have $\nabla \cdot \boldsymbol{v}_h \in \prob_{\polOrder}(\element)$. 
The boundary terms are also computable using the edge projection.
The VEM tuples for this example are summarised in \cref{tab: curl free vem}.
\begin{table}[htbp]
    \footnotesize
    \caption{VEM tuple summary for the curl free example. We use $I_{2 \times 2}$ to denote the identity matrix in $\R^2$.}\label{tab: curl free vem}
    \begin{center}
      \begin{tabular}{|c|c|c|c||c|c|} \hline
       $\valueBasisSet^{\element}$ & $\elementDofSet$ & $\cC^{\element}_0$ & $\gradBasisSet^{\element}$ & $\valueBasisSet^{s}$ & $\cC^s_0$ \\ 
       \hline
        $ \nabla \cM_{\polOrder+1}(\element)$ & 
        \cref{def: curl free dofs} & 
        \cref{eqn: curl free constraints}-\cref{eqn: extra curl free constraint} & 
        $\cM_{\polOrder}(\element)I_{2 \times 2}$ &
        $\cM_{\polOrder}(s) \boldsymbol{n} $ & 
        \cref{eqn: edge curl free constraints} \\
        \hline
      \end{tabular}
    \end{center}
\end{table}

\begin{remark}
    Similarly to \cref{rmk: trace of grad proj}, we can show ${\rm tr} (\gradProj v_h) \in \prob_{\polOrder}(\element)$. The projections are curl free due to the choice of basis sets $\valueBasisSet^{\element}:=\nabla \cM_{\polOrder+1}(\element)$ and ${\gradBasisSet^{\element} := \cM_{\polOrder}(\element)I_{2\times2}}$.
\end{remark}

\section{Extension to fourth-order problems}\label{sec: fourth order problems}
In this section we show how to generalise our framework to include fourth-order problems.
We restrict our extension of the generic framework to the $H^2$-conforming and nonconforming VEM spaces considered in \cite{antonietti_conforming_2018,brezzi_virtual_2013} and \cite{antonietti_fully_2018,10.1093/imanum/drab003,zhao_morley-type_2018}, respectively.
\changes{Additional nonconforming spaces which can be constructed using the concepts discussed in this paper were analysed and tested in \cite{10.1093/imanum/drab003}. These are especially suited for solving fourth-order perturbation problems.}


\subsection{Hessian projection}\label{sec: hessian projection}
Since we are now considering biharmonic-type problems, we need to introduce a further projection for the hessian term\corrections{; this will be denoted by} $\hessProj$. 
To this end, we extend the VEM tuple $\cV^{\element}$ to include an additional basis set, denoted by $\hessBasisSet^{\element}$ and assume in the following that $\hessBasisSet^{\element} \subset [\cM_{\polOrder-2}(\element)]^{2 \times 2}$.
\corrections{Therefore our revised VEM tuple is of the form}
\begin{align}\label{eqn: fourth order VEM tuple}
    \cV^{\element} = (\element, \valueBasisSet^{\element}, \elementDofSet, \cC^{\element}_0, (\cB^{\element}_i)_{i=1}^{2} ).
\end{align}
Similar to the gradient projection defined in \cref{def: gradient projection}, we define a fully computable hessian projection hierarchically using integration by parts once.

\begin{definition}[Hessian projection]\label{def: hessian projection}
    For a discrete function $v_h$ 
    we define the \emph{hessian projection} $\hessProj$ into the span of $\hessBasisSet^{\element}$ provided by the VEM tuple \eqref{eqn: fourth order VEM tuple}.
    \begin{multline}
        \int_{\element} \hessProj v_h : q \, \dx =
        - \sum_{i,j=1}^2 \int_{\element} [\gradProj v_h]_i \partial_j q_{ij} \, \dx
        \\
        + \sum_{s \in \cS_h(\element)} \int_{s} \left( \faceNormalProj v_h \sum_{i,j=1}^2 n_i n_j q_{ij}  + \Pi^s_1 v_h \sum_{i,j=1}^2 \tau_i n_j q_{ij} \right) \, \ds 
    \end{multline}
    for all $q \in \hessBasisSet^{\element}$. 
    Here, $n,\tau$ denote the unit normal and tangent vectors of $s$, respectively. 
\end{definition}
As you can see from \cref{def: hessian projection}, for the boundary terms we have used
two new projections, one for the tangential derivative $\Pi^s_1$ and one for the normal derivative $\faceNormalProj$ of $v_h$ on the edges. 
In 2D  we can use $\Pi^s_1 = \partial_s \faceProj$.
We will give details of the new normal derivative projection $\faceNormalProj$ for both of the examples considered in this section.





\subsection{Example: \texorpdfstring{$H^2$}{H2}-conforming VEM}\label{sec: H2 conforming vem}
The first example we consider is the $H^2$-conforming space and we take as $\elementDofSet$ the same dof set as in \cite{antonietti_conforming_2018}. 
First, as in \cite{antonietti_conforming_2018}, we assume $\polOrder\geq3$ but we will describe the lowest order case $\polOrder=2$, used for example in \cite{antonietti_$c^1$_2016}, at the end of the section.
\begin{definition}\label{def: H2 conforming dofs}
    For $\element \in \mesh$ the conforming dof set is given as follows.
    \begin{align}
        &\text{For $\polOrder \geq 2$: $h^{|p|}_{v} D^{p} v_h(v), \, |p| \leq 1,$ for any vertex $v$ of $\element$},
        \label{eqn: h2 conf vertex dofs}
        \\
        &\text{For $\polOrder \geq 4$: }
        \frac{1}{|s|} \int_s v_h m_{\polOrder-4} \, \ds \quad \forall m_{\polOrder-4} \in \cM_{\polOrder-4}(s), \ \text{ for any edge } s \in \cS_h(\element),
        \label{eqn: h2 conf edge dofs}
        \\
        &\text{For $\polOrder \geq 3$: }
        \int_{s} (\partial_n v_h) m_{\polOrder-3} \, \ds \quad \forall m_{\polOrder-3} \in \cM_{\polOrder-3}(s), \ \text{ for any edge } s \in \cS_h(\element),
        \label{eqn: h2 conf normal dofs}
        \\
        &\text{For $\polOrder \geq 4$: } 
        \frac{1}{|\element|} \int_{\element} v_h m_{\polOrder-4} \, \dx \quad \forall m_{\polOrder-4} \in \cM_{\polOrder-4}(\element).
    \end{align}
    Here, $h_v$ denotes a local length scale associated to the vertex v, e.g., an average of the diameters of all surrounding elements.
\end{definition}
The element value projection is defined exactly as before (\cref{def: value projection}) with constraints of the following form
\begin{align}\label{eqn: h2 conf value constraints}
    \cC(v_h) &:= \int_{\element} v_h m_{\polOrder-4} \, \dx, \quad m_{\polOrder-4} \in \cM_{\polOrder-4}(\element).
\end{align}

\begin{definition}[Edge value projection]\label{def: H2 conforming dofs edge value projection}
    We define the edge value projection as the unique solution $\faceProj v_h \in \text{Span}(\cM_{\polOrder}(s))$ of the following CLS problem. The least squares is made up of the following parts
    \begin{alignat*}{3}
        &\faceProj v_h (s^{\pm}) - v_h (s^{\pm}), &\qquad& 
        &D (\faceProj v_h (s^{\pm})) \cdot \tau - D v_h (s^{\pm}) \cdot \tau,
    \end{alignat*}
    where $s^{\pm}$ denotes the vertices attached to an edge $s$, subject to the constraints 
    \begin{align*}
        \int_s (\faceProj v_h - v_h) m^s_{\polOrder-4} \, \ds = 0, \quad m_{\polOrder-4}^s \in \cM_{\polOrder-4}(s).
    \end{align*}
\end{definition}
Even though we only have edge moments up to order $\polOrder-4$  \corrections{(which provide $\polOrder-3$ conditions)} we are able to define the edge value projection into $\prob_{\polOrder}(s)$.
\corrections{This is due to the dofs at the vertices. 
More specifically, the value and derivative dofs provide the extra conditions needed to determine a polynomial of degree $\polOrder$ along an edge $s$.}
Similarly, as shown in the next definition, despite only having edge normal dofs up to order $\polOrder-3$, we are able to define the edge normal projection into $\prob_{\polOrder-1}(s)$ due to the derivative values at the vertices.
\changes{In fact in both cases the CLS problem is equivalent to a square linear system of equations.}
\begin{definition}[Edge normal projection]\label{def: H2 conforming dofs edge normal projection}
    We define the edge normal projection as the unique solution $\faceNormalProj v_h \in \text{Span}(\cM_{\polOrder-1}(s))$ of the following CLS problem. The least squares is made up of the following parts
    \begin{alignat*}{3}
        &\faceNormalProj ( v_h (s^{\pm})) - D v_h (s^{\pm})\cdot n~,
    \end{alignat*}
    subject to the constraints
    \begin{align*}
        \int_s (\faceNormalProj v_h - \partial_n v_h) m_{\polOrder-3} \, \ds = 0, \quad m_{\polOrder-3}(s) \in \cM_{\polOrder-3}(s).
    \end{align*}
\end{definition}

The remaining parts of the VEM tuple $\cV^{\element}$ are summarised in \cref{tab: h2 conforming vem}.
\begin{table}[htbp]
    \footnotesize
    \caption{VEM tuple $\cV^{\element}$ summary for the $H^2$-conforming example (and $H^2$-nonconforming example). In particular, we take the hessian basis set to be the subset consisting of symmetric matrices of $[\cM_{r}(\element)]^{2 \times 2}$. \changes{We introduce $q,r \in \N$ for the gradient and hessian basis sets with $q \leq \polOrder$ and $r \leq \polOrder$.}}\label{tab: h2 conforming vem}
    \begin{center}
      \begin{tabular}{|c|c|c|} \hline
       $\valueBasisSet^{\element}$ & $\gradBasisSet^{\element}$ & $\hessBasisSet^{\element}$ \\ 
       \hline
        $\cM_{\polOrder}(\element)$ & 
        $[\cM_{q}(\element)]^2$ & 
        ${\rm sym}\big([\cM_{r}(\element)]^{2 \times 2}\big)$ 
        \\
        \hline
      \end{tabular}

    \end{center}
\end{table}
\changes{
Similar to \cref{rmk: stab free}, in the spirit of \cite{berrone2021lowest} we can take $q:=\polOrder$ and $r:= \polOrder-1$ in order to use a stabilised free fourth-order VEM space whilst $q:=\polOrder-1$ and $r:=\polOrder-2$ are the standard choices.}
Additionally, as stated in \cref{rmk: smaller basis sets}, for the $H^2$-spaces we can also use the smaller basis sets $\nabla\cM_{\polOrder}(\element)$ and $\nabla^2\cM_{\polOrder}(\element)$ which will work optimally for constant coefficient linear problems.

We conclude this example by describing how a minimal order $H^2$-conforming space can be constructed following \cite{antonietti_$c^1$_2016}:
the degrees of freedom are the same as before (given in \cref{def: H2 conforming dofs}).
Due to the fact that $\polOrder=2$, we have no inner moments and no edge moments but only values and derivatives at the vertices. 
As in the higher order case, this allows us to compute an \emph{edge normal projection} into $\prob_{1}(s)$ as detailed in \cref{def: H2 conforming dofs edge normal projection}.
Due to having a value and tangential derivative at both vertices of the
edge, the edge value projection needs to be computed into $\prob_{3}(s)$, i.e.,
into $\prob_{\polOrder+1}(s)$ instead of $\prob_{\polOrder}(s)$ as
described in \cref{def: H2 conforming dofs edge value projection}.

\corrections{For the element value projection,} using $\valueBasisSet^{\element} = \cM_{2}(\element)$ as in the
higher order case does not lead to a convergent method. Instead, our
numerical experiments indicate that the correct choice is the one given in
Table~\ref{tab: lowest order h2 conforming vem}.
\begin{table}[htbp]
    \footnotesize
    \caption{VEM tuple $\cV^{\element}$ summary for the lowest order $H^2$-conforming example.}
    \label{tab: lowest order h2 conforming vem}
    \begin{center}
      \begin{tabular}{|c|c|c|} \hline
       $\valueBasisSet^{\element}$ & $\gradBasisSet^{\element}$ & $\hessBasisSet^{\element}$ \\ 
       \hline
        $\cM_{3}(\element)$ & 
        $[\cM_{2}(\element)]^2$ & 
        ${\rm sym}\big([\cM_{1}(\element)]^{2 \times 2}\big)$ 
        \\
        \hline
      \end{tabular}
    \end{center}
\end{table}
An issue with the choice ${\valueBasisSet^{\element} = \cM_3(\element)}$ is that on triangles the space has nine degrees
of freedom (three per vertex) but the dimension of $\prob_{3}(\element)$
equals ten. Consequently, to uniquely define the value projection with our \refereecommentsvisible{constrained} least squares approach, we need to add a constraint.  In the absence of any inner moments,
we decided to introduce the following constraint in our implementation:
\begin{align*}
   \int_\element \Delta \valueProj v_h \, \dx = 
   \sum_{s \in \cS_h(\element)} \int_{s} \faceNormalProj v_h \, \ds.
\end{align*}
In future releases we will investigate adding more constraints of this form to define the value projection also for other spaces.


\subsection{Example: \texorpdfstring{$H^2$}{H2}-nonconforming VEM}\label{section: h2 nonconforming vem}
We now study the nonconforming space and provide the remaining quantities of the VEM tuple. The dof set $\elementDofSet$ is provided in \cref{def: H2 nonconforming dofs} and for this example we assume $\polOrder \geq 2$. 
\begin{definition}\label{def: H2 nonconforming dofs}
    For $\element \in \mesh$ the nonconforming dof set is given as follows.
    \begin{align}
        &\text{For $\polOrder \geq 2$: the values  of $v_h$ for any vertex $v$ of $\element$},
        \label{eqn: h2 nc vertex dofs}
        \\
        &\text{For $\polOrder \geq 3$: }
        \frac{1}{|s|} \int_s v_h m_{\polOrder-3} \, \ds \quad \forall m_{\polOrder-3} \in \cM_{\polOrder-3}(s), \ \text{ for any edge } s \in \cS_h(\element),
        \label{eqn: h2 nc edge value dofs}
        \\
        &\text{For $\polOrder \geq 2$: }
        \int_{s} (\partial_n v_h) m_{\polOrder-2} \, \ds \quad \forall m_{\polOrder-2} \in \cM_{\polOrder-2}(s), \ \text{ for any edge } s \in \cS_h(\element),
        \label{eqn: h2 nc normal edge dofs}
        \\
        &\text{For $\polOrder\geq4$: }
        \frac{1}{|\element|} \int_{\element} v_h m_{\polOrder-4} \, \dx \quad \forall m_{\polOrder-4} \in \cM_{\polOrder-4}(\element).
        \label{eqn: h2 nc inner dofs}
    \end{align}
\end{definition}
The element value projection is defined exactly as before (\cref{def: value projection}) with the same constraints as the previous conforming example \cref{eqn: h2 conf value constraints}.

Since we do not have derivative vertex dofs for this example, the edge value projection is also defined exactly as in \cref{def: value projection} as always with $\interfaceDofSet:=\Lambda^{\element,s}$,
where $\Lambda^{\element,s}$ contains the edge dofs in \cref{eqn: h2 nc edge value dofs} and the two vertex value dofs \cref{eqn: h2 nc vertex dofs}.
Note that we are able to define the edge value projection into $\text{Span}(\cM_{\polOrder-1}(s))=\prob_{\polOrder-1}(s)$ due to the additional vertex value dofs in $\Lambda^{\element,s}$, despite only having edge dofs up to order $\polOrder-3$. 
The constraints $\cC^{s}_0$ for the edge value projection are of the following form
\begin{align}
    \cC(v_h) &:= \int_{s} v_h m_{\polOrder-3} \, \ds, \quad m_{\polOrder-3} \in \cM_{\polOrder-3}(s).
    \label{eqn: h2 nc value edge constraints}
\end{align}

We define the edge normal projection into $\text{Span}(\cM_{\polOrder-2}(s))=\prob_{\polOrder-2}(s)$ with the definition given next, but note that this is enough for the hessian projection to be the exact $L^2$-projection of $\nabla^2 v_h$ as discussed in \cite{10.1093/imanum/drab003}.
\begin{definition}[Edge normal projection]
    We define the edge normal projection as the unique solution $\faceNormalProj v_h \in \text{Span}(\cM_{\polOrder-2}(s))$ of the following problem
    \begin{align*}
        \int_s (\faceNormalProj v_h - \partial_n v_h) m_{\polOrder-2} \, \ds = 0, \quad m_{\polOrder-2}(s) \in \cM_{\polOrder-2}(s).
    \end{align*}
\end{definition}

\corrections{The basis sets for the element value, gradient, and hessian projections are identical to those in the conforming example (\Cref{sec: H2 conforming vem}). 
Therefore, \cref{tab: h2 conforming vem} also describes the remaining VEM tuple quantities for this example.}

\ifthenelse{\boolean{thesiscorrections}}{\color{red}}{\color{black}}
\section{Extension to three dimensions}\label{sec: ext to 3D}
In this section we briefly describe how the framework can be extended to three space dimensions, \changes{under \Cref{assumption: dofs}}.
We first point out that the VEM tuple defined in \Cref{def: vem tuple} remains unchanged with the exception that $F \in \edgeSet$ now represents a face in 3D as opposed to an edge in 2D.
As before, for each $\element \in \mesh$, we build an element value projection $\valueProj$ and projections on each \refereecommentsvisible{two} dimensional interface (\Cref{def: value projection}) i.e. face projections $\edgeProj$ for each face $s \in \edgeSet$.
Using these projections, we then build a gradient projection $\gradProj$ analogously to the definition given in \Cref{def: gradient projection}.
We illustrate this extension by extending the example described in \Cref{sec: h1 VEM 2D} to three space dimensions.


\subsection{Example: \texorpdfstring{$H^1$}{H1}-conforming VEM in 3D}
For this example, the dofs are defined recursively using the $d=2$ dofs.
For each $\element \in \mesh$, we note that $\partial \element$ is a two dimensional face and so the VEM tuple for $s \in \edgeSet$ is given already in the 2D case in Example~\ref{2.1 dofs}. 
We take the degrees of freedom for the 3D space to be those in \Cref{def: h1 conforming dofs 3D}, provided next.

\begin{definition}\label{def: h1 conforming dofs 3D}
    For $\polOrder \geq 1$, we define the following dofs
    \begin{itemize}
        \item The dof set $\Lambda^{s}$ for $s \in \edgeSet$.
        \item For $\polOrder>1$, the moments of $v_h$ up to order $\polOrder-2$ inside the element $\element$
        \begin{align*}
            \frac{1}{|\element|} \int_{\element} v_h m_{\polOrder-2} \, \dx \quad &\forall m_{\polOrder-2} \in \cM_{\polOrder-2}(\element).
        \end{align*}
    \end{itemize}
\end{definition}

We also note that the recursive construction of this 3D space means that the face projection in 3D is identical to the element value projection in 2D. 
Therefore, there is no difference between the constraints when $d=2,3$ for both value projections, other than whether $s \in \cS_h(\element)$ represents an edge or a face. 
Therefore, the constraints $\cC^s_0$ for the face projections are given by \eqref{eqn: h1 conforming edge constraints} and the constraints $\cC^{\element}_0$ for the element value projection are given by \eqref{eqn: h1 conf value constraints}.
The VEM tuples for this 3D example are described in \Cref{tab: h1 conforming vem tuple 3D}.
\begin{table}[htbp]
    \footnotesize
    \caption{VEM tuple summary for the $H^1$-conforming VEM in three dimensions.}\label{tab: h1 conforming vem tuple 3D}
    \begin{center}
      \begin{tabular}{|c|c|c|c||c|c|} \hline
       $\valueBasisSet^{\element}$ & $\elementDofSet$ & $\cC^{\element}_0$ & $\gradBasisSet^{\element}$ & $\valueBasisSet^{s}$ & $\cC^s_0$ \\ 
       \hline
        $\cM_{\polOrder}(\element)$ &
        \cref{def: h1 conforming dofs 3D} &
        \cref{eqn: h1 conf value constraints} &
        $[\cM_{\polOrder-1}(\element)]^2$ &
        $\cM_{\polOrder}(s)$ &
        \cref{eqn: h1 conforming edge constraints}
        \\
        \hline
      \end{tabular}
    \end{center}
\end{table}
\color{black}

\section{Implementation details}\label{sec: implementation details}
In this section we give an overview of our approach to implementing
VEM spaces within an existing finite element software framework.
\changes{The concepts for implementing general VEM spaces described here should easily carry over to other FE packages, e.g., \cite{alnaes2015fenics,bangerth2007deal,hecht2012new,prud2012feel,rathgeber2016firedrake} since it is minimally invasive.}
\corrections{The implementation has been made available in} the \textsc{Dune-Fem} \cite{dedner2010generic} module which is part of the Distributed Unified Numerics Environment DUNE \cite{dunegridpaperII}. 
So far, we have a proof of concept implementation \corrections{of all} spaces described in the previous section \corrections{only} in two space dimensions. 
We believe that most of the required steps readily carry over to higher space dimensions.
We can only give a very short overview here, \corrections{assuming} that the reader is to a certain extent familiar with the implementation of finite element methods on unstructured grids.
\changes{In most cases, the implementation of general FEM spaces is based on the triple
$(\element,\cB^\element,\Lambda^\element)$ consisting of:
\begin{itemize}
    \item An element $\element$ in an appropriate tessellation $\mesh$ of the domain $\Omega$ in which the problem is posed.
    \item A basis $\cB^\element := ( b_{\alpha}^\element )_{\alpha=1}^{N^{\element}}$ of the finite element space, where $N^{\element}$ denotes the dimension of the FEM space
          \changes{(more commonly the space is used directly as a component of the triple).}
    \item A set of functionals (degrees of freedom) $\elementDofSet := ( \lambda_j^\element )_{j=1}^{N^{\element}}$.
\end{itemize}
The degrees of freedom are assumed to be unisolvent, \refereecommentsvisible{i.e. if $\Lambda^{\element} (v_h) = 0$ then $v_h \equiv 0$.}
Consequently, the \refereecommentsvisible{local} matrix $\tilde A^\element$ of size $(N^{\element} \times N^{\element})$ defined by \refereecommentsvisible{evaluating the degrees of freedom at the finite element basis, i.e.}
$\tilde A^{\element}_{i, \alpha} = \lambda_i^\element (b_{\alpha}^\element)$ is regular.
The inverse of $\tilde A^{\element}$ can be used to construct the nodal basis functions $\Phi^\element = ( \phi_j^\element )_{j=1}^{N^{\element}}$ on each element $E$, which satisfy
$\lambda_i^\element (\phi_j^\element) = \delta_{ij}$,
and are given by \cite{10.1007/978-3-642-28589-9_1}
\begin{align}\label{eqn: nodal basis defn}
    \Phi^\element(x) = ( \tilde A^\element)^{-1}\cB^\element(x).
\end{align}
Code to evaluate the nodal basis $\Phi^\element$ and their derivatives for a given FEM triple forms the fundamental building block of most finite element packages. In the following we will demonstrate that the implementation of VEM spaces can be done following the same concepts using our VEM tuples which are a direct extension of the FEM triple.
Note that the main difference between the finite element triple and the virtual element tuple is that the basis set $\valueBasisSet^{\element}$ is not a basis for the virtual element space.
In particular, we can have $\text{dim} (\valueBasisSet^{\element}) < N^{\element}$ where $N^{\element}$ denotes the dimension of the discrete space.
This means that the matrix $\elementDofSet(\cB^\element)$ is no longer square and the nodal basis can not be directly computed as in the FE setting.}

\subsection{General structure of assembly code}
Most finite element implementations will in some form or another contain the following.
\begin{itemize}
\item A way to iterate over the tessellation $\mesh$, e.g., the triangles of a simplex grid. 
\item For each element $\element \in \mesh$ and given order $p$, some form of quadrature rule 
      \[ 
        Q^E_p(f)=\sum_{q=1}^{N_Q}\omega^E_q f(x^E_q),
      \]
      where $\omega^E_q$ are the weights and $x^E_q\in E$ are the quadrature points.
\item The evaluation of the nodal basis $\Phi^E=(\phi^E_k)_{k=1}^{N^E}$
      of the local finite element space and its derivatives. 
      The definition \corrections{of $\Phi^E$} is given in \cref{eqn: nodal basis defn} and satisfies $\Lambda^{\element}(\Phi^{\element})=I$.
\item A \emph{local to global} dof mapper $\mu^E=(\mu^E_k)_{k=1}^{N^{\element}}$.
      This takes the numbering
      of the local degrees of freedom and converts them into some global index
      so that shared dofs on faces, edges, and vertices correspond to the
      same entry in some global storage structure for the degrees of
      freedom.
\end{itemize}
Using the above, the assembly of for example a functional of the form 
${b(v)=(f,v)_{\Omega}}$ is \corrections{showcased} in \cref{AlgAssembly}.
\begin{algorithm}
  \caption{Assembly of a functional of the form $b(v)=(f,v)_{\Omega}$}
  \label{AlgAssembly}
  \begin{algorithmic}[1]
  \STATE Set the vector $b$ to zero
  \FOR{each element $\element \in \mesh $}
  \STATE Compute the local contribution $b_k^E = Q^E_p(f\phi^E_k)$ for $k=1,\dots,N^E$
  \STATE Scatter the local contributions into the global vector $b$, $b_{\mu^E_k} = b_{\mu^E_k} + b_k^E$
  \ENDFOR
  \end{algorithmic} 
\end{algorithm}

\textsc{Dune} can directly work on polygonal grids but for quadrature we need
to subdivide each polygon into  \corrections{standard domains such as} triangles.
For our implementation of the VEM spaces we therefore decided to directly work
with such a \corrections{sub-triangulation}.
So the iteration over the elements in \cref{AlgAssembly} would be over the sub-triangulation denoted by $T_h$. 
For each triangle $T\in T_h$ we have access to the unique polygon $E_T$ with $T\subset E_T$. 
This can be seen as a colouring or \emph{material property} of each triangle. 

Now, the above algorithm can be used unchanged with a VEM space: the dof mapper $\mu^T$ for a triangle $T$ is defined to be $\mu^{E_T}$ and instead of
$\phi^T_k$ one uses $\valueProj\phi^{E_T}_k$. 
\corrections{For instance,} if the existing assembly code is based on a function to evaluate all basis functions $(\phi^T_{k}(x))_k$ at a given point $x\in T$, one needs to provide a function to compute
$(\valueProj\phi^{E_T}_k(x))_k$ to use a VEM space.
In the same way the assembly of the stiffness matrix requires a function evaluating the gradients which in the FEM case would return $(\nabla\phi^T_k(x))_k$ \corrections{while} in the VEM case it would return
$(\gradProj\phi^{E_T}_k(x))_k$.

\subsection{Nodal basis functions and projection operators}
As the above demonstrates, very little change is required to the assembly step of a finite element package \corrections{to include our VEM implementation}.  
The main new part required is code to compute $(\valueProj\phi^E_k(x))_k$ and $(\gradProj\phi^E_k(x))_k$
for a given polygon $\element \in \mesh$ and $x\in \element$. 
We will give a brief summary for $\valueProj$; the projections for the higher derivatives, e.g., $\gradProj$, are even simpler and identical for all spaces (see \cref{def: gradient projection} and \cref{def: hessian projection}).

The range space for $\valueProj$ is spanned by
$\valueBasisSet^{\element}=( m^E_\alpha )_\alpha$ and is a subspace of the polynomial space $\prob_{\polOrder}(\element)$ spanned
by the basis set $\cM_{\polOrder}(\element)$.
In our implementation we construct a \emph{minimal area bounding box} for each polygon $\element$ and use scaled tensor product Legendre basis functions over this rectangle to define $\cM_{\polOrder}(\element)$.  To increase stability of the system matrices we also provide the option to construct an orthonormalised version of these monomials for each $\element$ based on the $L^2$-scalar product over $\element$. 
\corrections{More information and discussion on how the choice of basis affects the resulting system can be found in \cite{mascotto_ill-conditioning_2018}.}

Next we construct $\valueProj\phi^{\element}_k$ in the span of
$\valueBasisSet^{\element}$.
To do this, we construct the matrix
$\boldsymbol{\Pi}^{\element}_0$ which provides the coefficients for each basis function $k=1,\dots,N^{\element}$, where
\begin{align*}
  \valueProj \phi^E_k = \sum_{\alpha=1}^{\mathcal{N}^{\element}_0} (\boldsymbol{\Pi}^{\element}_0)_{\alpha, k} m_{\alpha}^{\element}
\end{align*}
so that the $k$-th column of $\boldsymbol{\Pi}^{\element}_0$ contains the coefficients of the polynomial $\valueProj \phi^{\element}_k$ in the basis $\valueBasisSet^{\element} = (m^{\element}_\alpha)_{\alpha}$.
We denote with $\mathcal{N}^{\element}_0$ the size of the basis set $\valueBasisSet^{\element}$.
The matrix $\boldsymbol{\Pi}^{\element}_0$ can be computed based on the
\emph{\refereecommentsvisible{constrained} least squares problem} from \cref{def: value projection}.
\begin{align*}
    \sum_{i=1}^{\NoLocalDofs} (\lambda_i (\valueProj \phi^{\element}_k) - \lambda_i(\phi^{\element}_k))^2 
    &= 
    \sum_{i=1}^{\NoLocalDofs} \left( 
      \sum_{\alpha=1}^{\mathcal{N}^{\element}_0} (\boldsymbol{\Pi}^{\element}_0)_{\alpha,k}  \lambda_i ( m^{\element}_{\alpha} ) - \lambda_i(\phi^{\element}_k) \right)^2 
    \\
    &= \| \tilde A^{\element} (\boldsymbol{\Pi}^{\element}_0)_{k}  - \lambda_i(\phi^{\element}_k) \|^2
\end{align*}
where the matrix $\tilde A^{\element}$ is the same basis transformation matrix used in the FEM setting and
the vector $(\boldsymbol{\Pi}^{\element}_0)_{k}$ denotes the $k$-th column of $\boldsymbol{\Pi}^{\element}_0$.
Assuming that the nodal basis of the local VEM space satisfies $\lambda_i(\phi_j^\element) = \delta_{ij}$ for each $\lambda_i \in \Lambda^{\element}$, then it is clear that $\lambda_i(\phi^{\element}_k)$ reduces to the $k$-th unit vector.

It remains to describe the set up of the constraint system,  \corrections{that is} the remaining part of \cref{def: value projection}. 
Recall that these are described as follows
\begin{align*}
  \cC ( \valueProj \phi^{\element}_k - \phi^{\element}_k ) =
  \cC \left( \sum_{\alpha=1}^{\mathcal{N}^{\element}_0} (\boldsymbol{\Pi}^{\element}_0)_{\alpha, k} m_{\alpha}^{\element} - \phi^{\element}_k \right) =
  \sum_{\alpha=1}^{\mathcal{N}^{\element}_0} (\boldsymbol{\Pi}^{\element}_0)_{\alpha, k} \cC(m_{\alpha}^{\element}) - \cC(\phi^{\element}_k) = 0.
\end{align*}

\ifthenelse{\boolean{thesis}}
{
  For the CLS problem with $A \in \R^{n\times m}$, $B \in \R^{p \times n}$. 
  The CLS problem is as follows: find $x \in \R^n$ which solves 
  \begin{align}
    {\rm min}_{x} \| Ax -b\|_{2} \quad \text{ such that } \quad Bx=d
  \end{align}
  which is uniquely solvable if 
  \begin{align*}
    {\rm rank}(B) = p \quad \text{ and } \quad {\rm rank} \left( \begin{array}{c}  A 
      \\
      B 
    \end{array}
    \right) = n 
  \end{align*} 
}
{
}

\begin{remark}\label{rmk: CLS solvability}
  We note that all of the constraint matrices $\tilde C$ for the examples in \Cref{sec: examples} and \Cref{sec: fourth order problems} are of the form 
  \begin{align}\label{eqn: mass matrix}
    \tilde C_{\beta \alpha} = \int_{\element} m_{\beta}^{\element} m_{\alpha}^{\element} \, \dx 
  \end{align}
  for $m_{\alpha}^{\element} \in \valueBasisSet^{\element}$ and for $m_{\beta}^{\element} \in \cB^{\element}_*$ where $\cB^{\element}_* \subset \valueBasisSet^{\element}$.
  \corrections{For instance,} in Example~\ref{constraint sets} we take $\cB^{\element}_* := \cM_{\polOrder-2}(\element)$.
  Since the constraint matrices in \cref{eqn: mass matrix} are truncated mass matrices of size $(\mathcal{N}^{\element}_{*} \times \mathcal{N}^{\element}_0)$ where $\mathcal{N}^{\element}_{*} := {\rm dim}(\cB^{\element}_*)$ they have full rank.
  The CLS problem is uniquely solvable
  if and only if (e.g. \cite{bjorck1996numerical})
  \begin{align*}
    {\rm rank}(\tilde C) = \mathcal{N}^{\element}_{*} \quad \text{ and } \quad  
    {\rm rank} \left( \begin{array}{c} \tilde A^{\element} 
      \\
      \tilde C
    \end{array}
    \right) =  \mathcal{N}^{\element}_0.
  \end{align*}
  The second condition is equivalent to ${\rm Null}(\tilde A^{\element}) \cap {\rm Null} (\tilde C) = \{ 0 \}$ shown as in e.g. \cite{bjorck1996numerical} but it is clear that ${\rm ker}(\tilde A^{\element}) = \{ 0 \}$ due to unisolvency of the dofs and using that $\prob_{\polOrder}(\element)$ is a subset of the local VEM space (\cref{assumption: dofs}).
  \ifthenelse{\boolean{thesis}}{
    If this was not the case, there would exist $x \in {\rm ker}(\tilde A^{\element})$ with $x \neq 0$ such that $\tilde A^{\element} x=0$ i.e. for each $i=1,\dots,N^{\element}$,
  \begin{align*}
    0 = \sum_{\alpha=1}^{\mathcal{N}^{\element}_0} \tilde A^{\element}_{i\alpha} x_\alpha = \sum_{\alpha=1}^{\mathcal{N}^{\element}_0} \lambda_i ( m^{\element}_{\alpha}) x_\alpha = \lambda_i \left( \sum_{\alpha=1}^{\mathcal{N}^{\element}_0} x_\alpha m_\alpha^{\element} \right)
  \end{align*}
  but since the polynomial space $\prob_{\polOrder}(\element)$ is a subset of the local VEM space, this contradicts the unisolvency assumption of the degrees of freedom. 
  }{}
\end{remark}

\subsection{Usage in \textsc{Dune-Fem}}
We provide Python bindings for our package that make use of Unified Form Language (UFL) \cite{alnaes_unified_2012} for the problem formulation. 
Therefore a problem of the form
$$ \int_\Omega D(x,u)\nabla u\cdot\nabla v \, \dx + \int_{\Omega} m(x,u)v \, \dx = 0 $$
using a standard Lagrange space on an unstructured grid is shown in Listing~\ref{AlgLagrange}.
\begin{lstlisting}[label = AlgLagrange, caption = Python script to set up a Lagrange \corrections{space} using \textsc{Dune-Fem}, language = python]
from ufl import TrialFunction,TestFunction,SpatialCoordinate,grad,dot,dx
import dune.fem, dune.alugrid
gridView = dune.alugrid.aluConformGrid( gridDict )
spc = dune.fem.space.lagrange(gridView, order=k)
u,v,x = TrialFunction(spc),TestFunction(spc),SpatialCoordinate(spc)
D, m = 1+u*u, 2*u + cos(u) - cos(dot(x,x)) # example non-linearity
eqn = ( D*dot(grad(u), grad(v)) + m*v ) * dx == 0
scheme = dune.fem.scheme.galerkin( eqn )
uh = space.interpolate(0,name="solution")
scheme.solve(target=uh)
\end{lstlisting}
The grid is constructed using a Python dictionary {\tt gridDict}
containing the points and connectivity for the tessellation of
the domain $\Omega$.
To use a VEM space on a polygonal grid, both the grid and space construction part
needs to be adapted as shown in Listing~\ref{AlgVEMSpace}.
\begin{lstlisting}[label=AlgVEMSpace, caption = Modifications to the Python script in order to set up a VEM space using \textsc{Dune-Vem}, language = python]
import dune.vem
gridView = dune.vem.polyGrid( gridDict )
spc = dune.vem.vemSpace( gridView, order=k, testSpaces=[0,k-2,k-2])
\end{lstlisting}

In the {\tt gridDict} the {\tt simplices} key is replaced by {\tt polygons}.
The {\tt testSpaces} argument in the space constructor can be used to set the vertex, edge, and
inner moments to use for the degrees of freedom. So the above defines a
conforming VEM space (see \Cref{sec: h1 VEM 2D}).
\examplechanges{Using {\tt testSpaces=[-1,k-1,k-3]} results in the simplest serendipity version of this space (see \cref{rem: serendipity VEM}).}
The nonconforming space, for example, is constructed if {\tt testSpaces=[-1,k-1,k-2]} was used.
The $H^2$-nonconforming space (see the example in \Cref{section: h2 nonconforming vem}) requires
{\tt testSpaces=[-1,[k-3,k-2],k-4]} where the second argument now defines
the value moments and the edge normal moments to use.
This is the same notation introduced as the dof tuple in \cite{10.1093/imanum/drab003}.
\changes{The default range spaces for the value, gradient, and hessian
projections are polynomial spaces of order $k,k-1,k-2$, respectively. 
To use different values, the {\tt
order} parameter can be used by passing in a list instead, e.g., {\tt order=[k,k,k-1]} provides the projection operators used in the example for the
non-stabilised methods in \Cref{sec: testStab}}.

Other available spaces at the time of writing are {\tt divFreeSpace} (see \Cref{sec: div free}) and {\tt curlFreeSpace}
(see \Cref{sec: curl free}).
Other spaces, e.g., $H(div)$ and $H(curl)$-conforming spaces will be added in a later release.

The final required change concerns the stabilisation. 
After assembly we need to add the stabilisation term given by the matrix $S^E$ and some scaling $ ( \bar{D} + \bar{m}h^2 )S^{\element},$ where $S^{\element}$ is given by the \refereecommentsvisible{\emph{dofi-dofi}} stabilisation \cite{beirao_da_veiga_basic_2013} and is independent of the problem.

To allow for some flexibility, especially in the nonlinear setting, the
{\tt scheme} constructor takes two additional arguments used for
$\bar{D},\bar{m}$, respectively. 
For the above problem one can simply use
the same UFL expressions used to define the bilinear form as shown in Listing~\ref{AlgVemScheme}.
\begin{lstlisting}[label=AlgVemScheme, caption = Additional arguments in the VEM scheme constructor to handle the stabilisation term, language = python]
dune.vem.vemScheme(eqn, gradStabilization=D, massStabilization=m)
\end{lstlisting}

\begin{remark}\label{rmk: dune-vem installation}
  The package can be simply downloaded through the Python Package Index
  (PyPI) using: 
  {\tt pip install dune-vem }
  and a number of examples can be found in the \textsc{Dune-Fem} tutorial \cite{dednerPythonBindings2020}.
\end{remark}

\section{Numerics}\label{sec: numerics}
In the following we \corrections{report} results for a Laplace problem, \changes{an investigation of the non-stabilised methods,} followed by two time dependent problems using the curl free, second-order, divergence free, and fourth-order spaces, respectively. 
The aim of this section is not a detailed discussion of convergence rates and performance
of the virtual element methods. 
This has been done in the literature already referenced
in the example section, \Cref{sec: examples}. 
For discussions on fourth-order problems using the presented code see~\cite{dedner2021higher,10.1093/imanum/drab003}.
We will therefore only present some results to demonstrate
the wide range of problems that can be tackled with the presented VEM
spaces. 
The code to run all examples in this section as well as detailed instructions on how to run the scripts can be found at the \texttt{dune-vem-paper} repository\footnote{\href{https://gitlab.dune-project.org/dune-fem/dune-vem-paper}{https://gitlab.dune-project.org/dune-fem/dune-vem-paper}}.

\begin{remark}
  \changes{
  In our examples we focus on problems with Dirichlet or \refereecommentsvisible{homogeneous} Neumann boundary conditions. 
  To add \refereecommentsvisible{non homogeneous} Neumann or even Robin boundary conditions, terms involving integrals over the domain boundary are
  required. 
  For example, consider the local form
  \begin{align*}
    a^\element(v,w) =
      \int_{\element} D(v,\nabla v) \cdot \nabla w \, \dx + \int_{\element} m(v,\nabla v) w \, \dx
      + \sum_{s \in \cS_h(\element)\cap\partial\Omega} \int_{s} (\alpha v-g_R)w \, \ds.
  \end{align*}
  To implement the discrete version of the boundary term requires the edge projection, 
  $$ \sum_{s \in \cS_h(\element)\cap\partial\Omega} \int_{s} \Big(\alpha \faceProj v-g_R\Big) \faceProj w \, \ds. $$
  \refereecommentsvisible{Note that we do not use the edge interpolation here since even though these coincide in most cases (in two space dimensions) they do not in the $H^2$-nonconforming case. 
  An investigation of another fourth-order VEM space exploring this property can be found in \cite{10.1093/imanum/drab003}.}
  Our numerical experiments indicate that using the element projections
  instead does not lead to a convergent scheme.  
  To the best of our knowledge, numerical analysis of these problems is not yet available.}

  \changes{We note that at the time of writing, we are restricted to \refereecommentsvisible{homogeneous} Dirichlet and Neumann boundary conditions for the fourth-order problems. We hope to extend this in further releases.}
\end{remark}

\subsection{Laplace problem: primal and mixed form}

\begin{figure}[!ht]\label{fig: Laplace}
  \resizebox{\textwidth}{!}{
\centering
\includegraphics[width=0.19\textwidth]{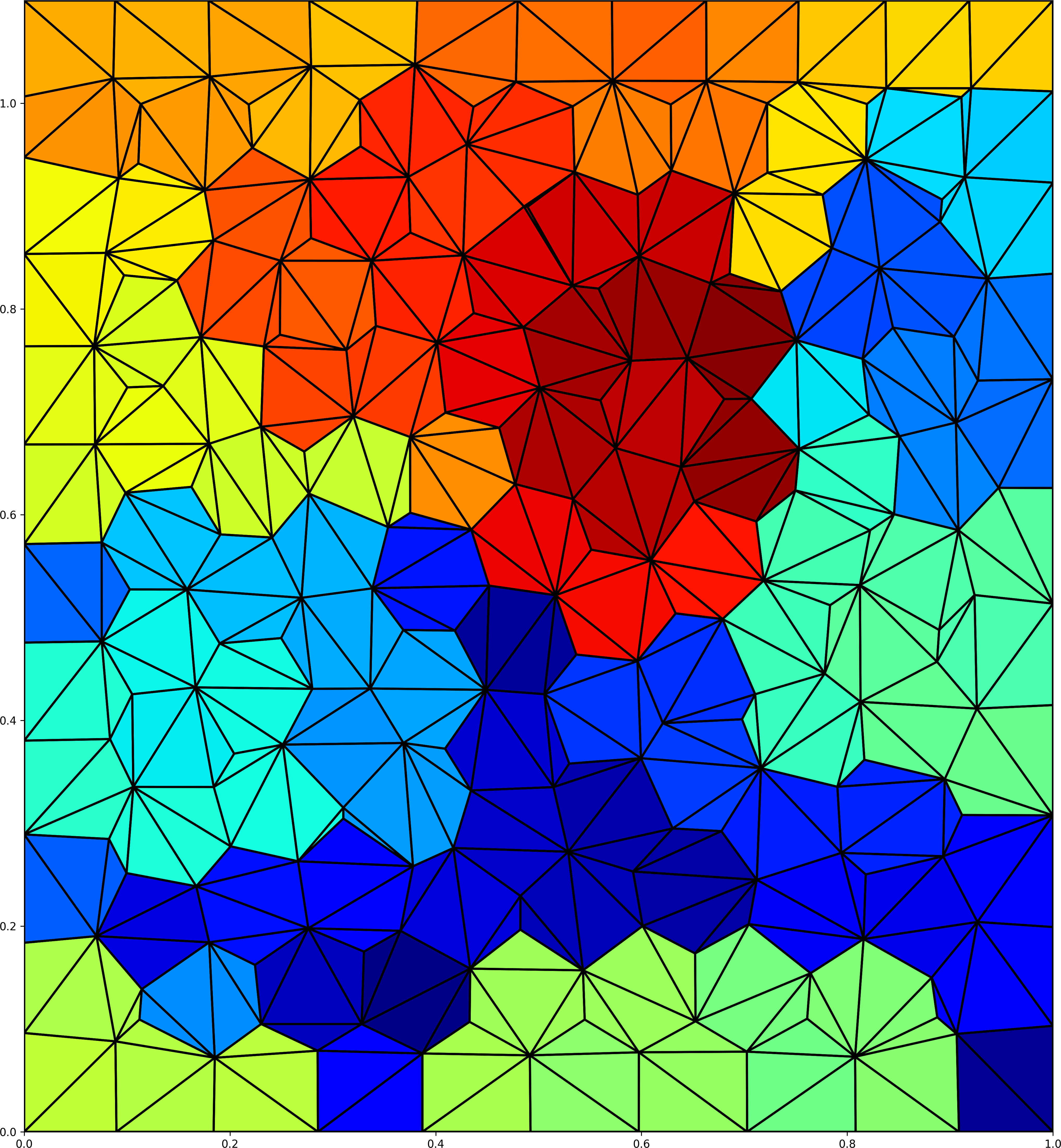} \hfill
\includegraphics[width=0.38\textwidth]{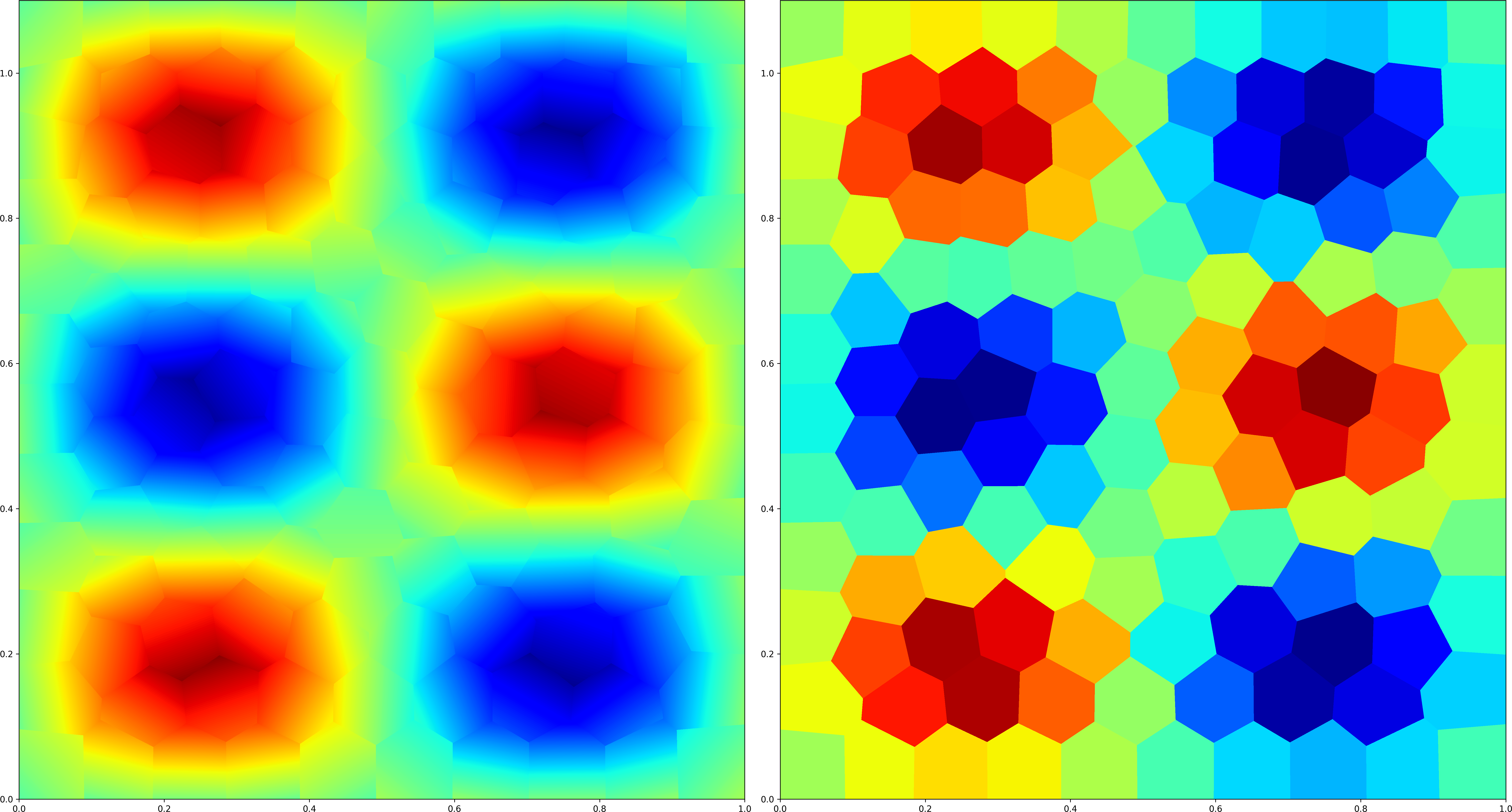} \hfill
\includegraphics[width=0.38\textwidth]{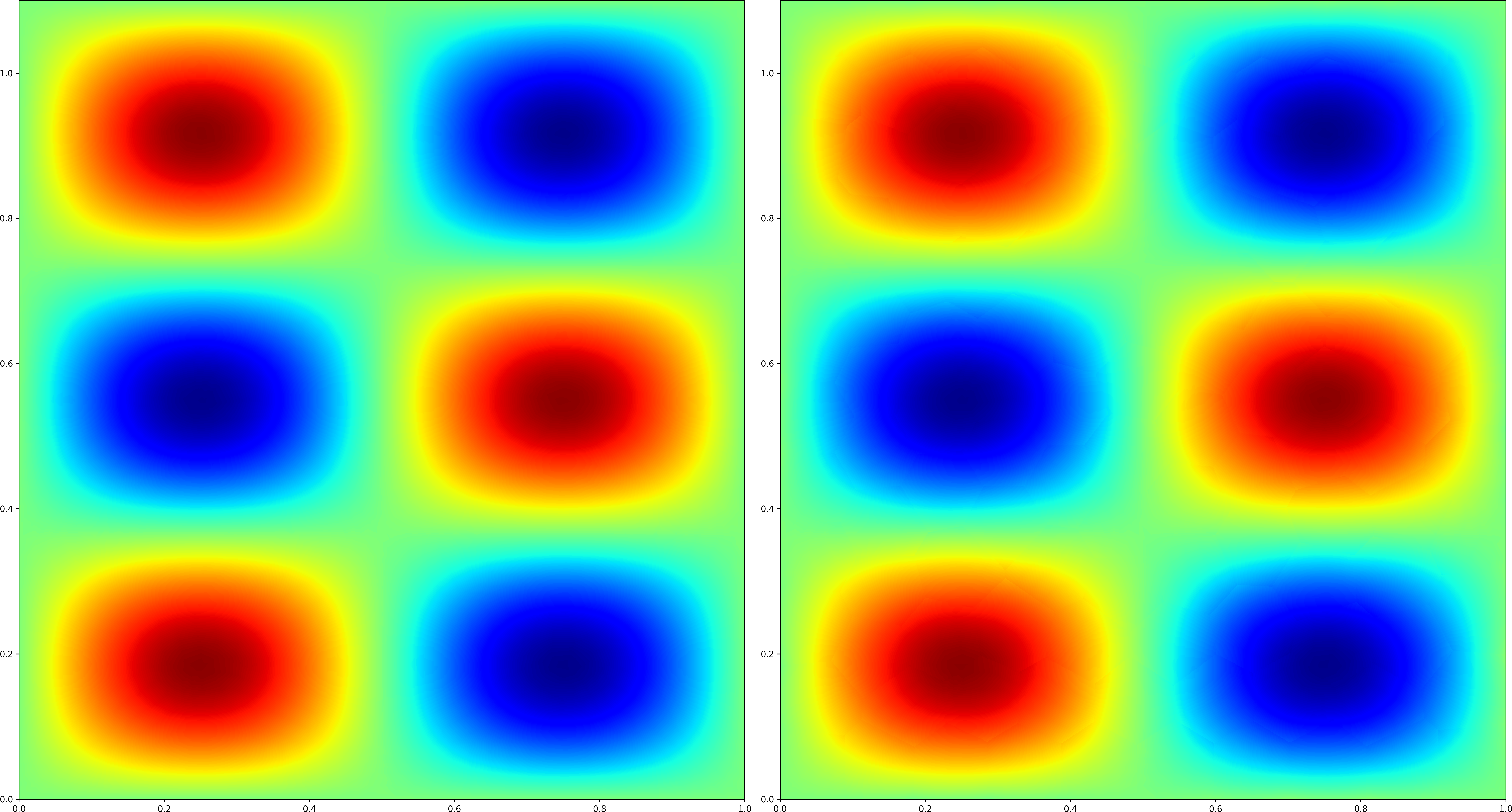}
  }
\caption{\changes{
Solutions to the Laplace problem using two different orders for the
approximation (left to right) using a Voronoi grid.
The left figure shows the grid together with the sub-triangulation used for our
agglomeration approach.
The next two figures show the solution using a conforming space
of order $1$ and the mixed approach with a DG and curl free space of order $0$, respectively. 
The right two figures show again the
primal and mixed solutions but with $\polOrder=2$ and $\polOrder=1$, respectively.}}
\end{figure}

In our first example we solve a simple Laplace problem with Dirichlet
boundary conditions on the square $[0,L_x]\times[0,L_y]$
subdivided into Voronoi cells.
We use two different approaches to solve the problem. Firstly, we use a standard primal formulation \corrections{as follows: find $u \in H^1_0(\Omega)$ such that} 
$$ a(u,v) := (\nabla u, \nabla v) = (f,v) \quad \forall v \in H^1_0(\Omega),$$ 
\corrections{using the conforming VEM space (\Cref{sec: h1 VEM 2D}).}
Secondly, we use a mixed formulation:
\corrections{find $(\sigma,u) \in H({\rm div},\Omega) \times L^2(\Omega)$ such that}
\changes{
\begin{align*}
   \int_\Omega\sigma\cdot\tau + u{\rm div}\tau \, \dx = 0~, \quad
   \int_\Omega ({\rm div}\sigma+f) v \, \dx = 0~,\quad
  \forall (\tau,v)\in H({\rm div},\Omega) \times L^2(\Omega)~.
\end{align*}
We discretise using a discontinuous Galerkin (DG) space for $u$
and the curl free space for the flux $\sigma$.}
The forcing $f$ is chosen so that the exact solution is given by
$\sin(2\pi x/L_x)\sin(3\pi y/L_y)$ and we use $L_x=1,L_y=1.1$.
Results are shown in \cref{fig: Laplace}.

\subsection{Non-stabilised methods}\label{sec: testStab}

In this next example we solve both a linear second-order and fourth-order
problem using the conforming spaces (\Cref{sec: h1 VEM 2D} and
\Cref{sec: H2 conforming vem}, respectively):
find $u\in H^p_0$ such that
\begin{align*}
  a(u,v) := (\kappa(\boldsymbol{x}) \nabla^p u, \nabla^p v) = (f,v) \quad \forall v \in H^p_0,
\end{align*}
with $p=1$ or $p=2$ and
$\kappa(\boldsymbol{x}) = \frac{10}{0.01+\boldsymbol{x}\cdot\boldsymbol{x}}$.
We set the forcing $f$ so that the exact solution is given by
$u(x,y) = (\sin(2\pi x)\sin(2\pi y))^2$ on the domain $\Omega=[0,1]^2$, and we use a Cartesian grid in both cases.

\begin{figure}[!ht]\label{fig: testStab}
  \begin{center}
    \includegraphics[width=0.49\textwidth]{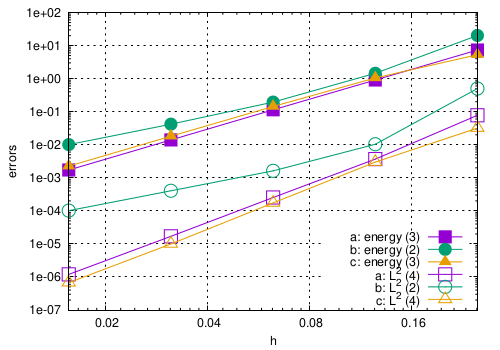}
    \includegraphics[width=0.49\textwidth]{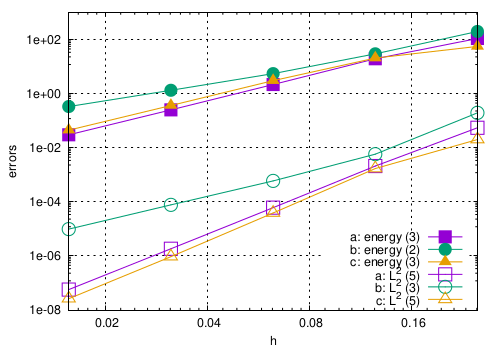}
  \end{center}
    \caption{Energy and $L^2$-norm errors for different grid resolutions using the stabilised and non-stabilised
    versions of the VEM discretisation. Left figure shows results for a second-order problem
    with the $H^1$-conforming VEM with order $\polOrder=3$ and the right figure shows results for a fourth-order problem with the $H^2$-conforming VEM with order $\polOrder=4$. Recall that (a) is the original space with stabilisation, (b) the same space
    without, and (c) is without stabilisation but with the extended range
    spaces for the gradient and hessian projections.}
\end{figure}

We test three different variants of the discretisation: (a) the standard
version of the spaces with the \refereecommentsvisible{\emph{dofi-dofi}} stabilisation \cite{beirao_da_veiga_basic_2013},
(b) the standard version of the spaces without any stabilisation, and finally (c)
a version of the spaces where the gradient and hessian projections are computed into 
$[\cM_{\polOrder}(\element)]^2$ and $[\cM_{\polOrder-1}(\element)]^{2\times 2}$
(instead of $\polOrder-1$ and $\polOrder-2$ as in the standard case).
At the time of writing there is no proof that this will lead to a stable
scheme but results shown in \cref{fig: testStab} indicate that the
non-stabilised method with the extended range spaces for the projections is
comparable to the original approach outperforming the non-stabilised method
used with the original space which leads to a suboptimal convergence rate.

\subsection{Incompressible flow around a cylinder}
\begin{figure}[!ht]\label{fig: navier stokes}
  \begin{center}
    \includegraphics[width=0.49\textwidth]{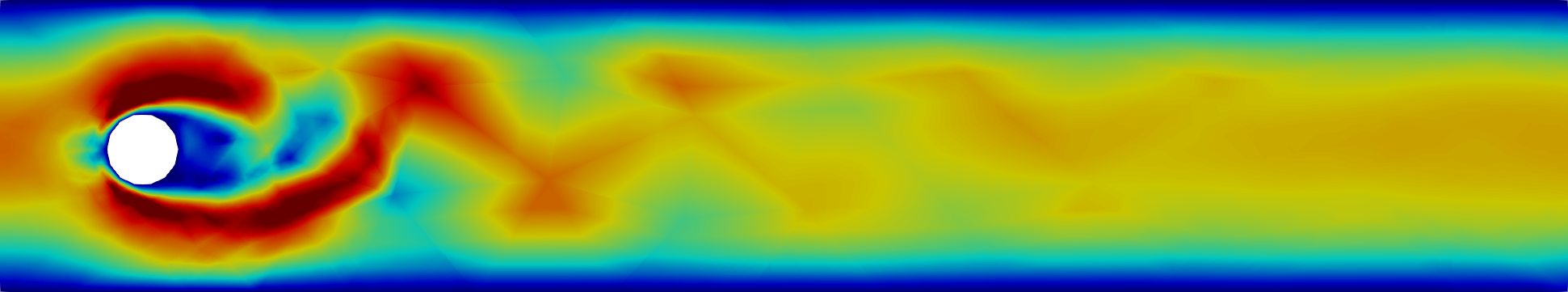}
    \includegraphics[width=0.49\textwidth]{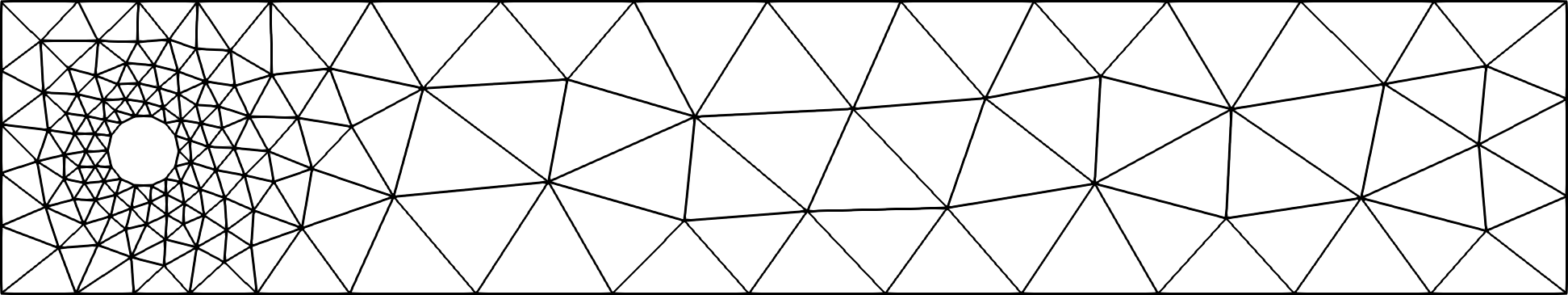} \\[0.2em]
    \includegraphics[width=0.49\textwidth]{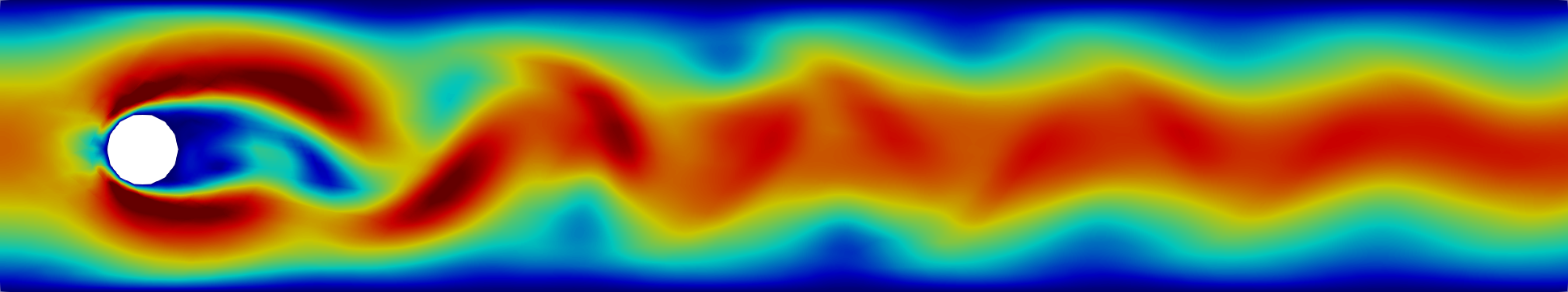}
    \includegraphics[width=0.49\textwidth]{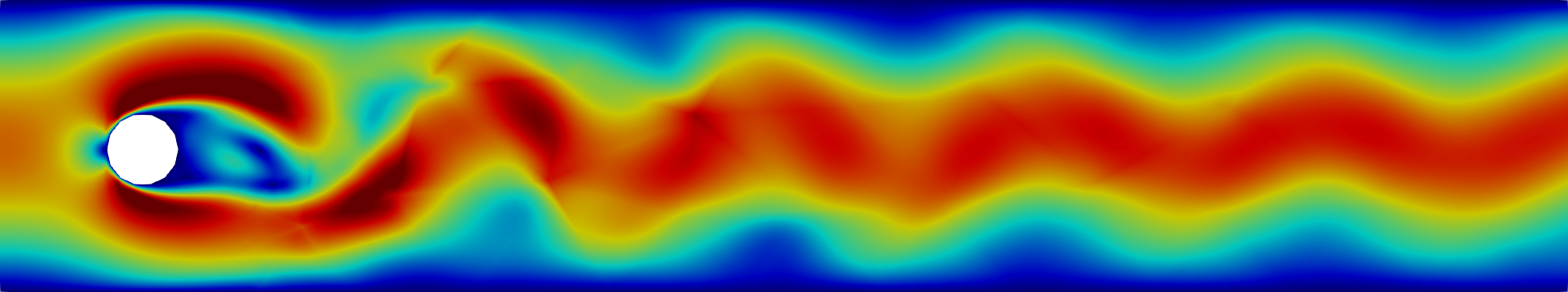} \\[0.2em]
    \includegraphics[width=0.49\textwidth]{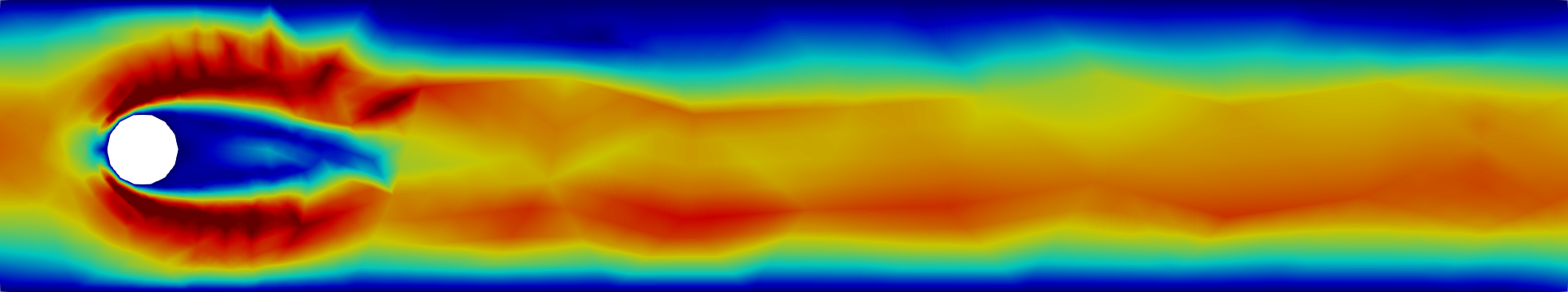}
    \includegraphics[width=0.49\textwidth]{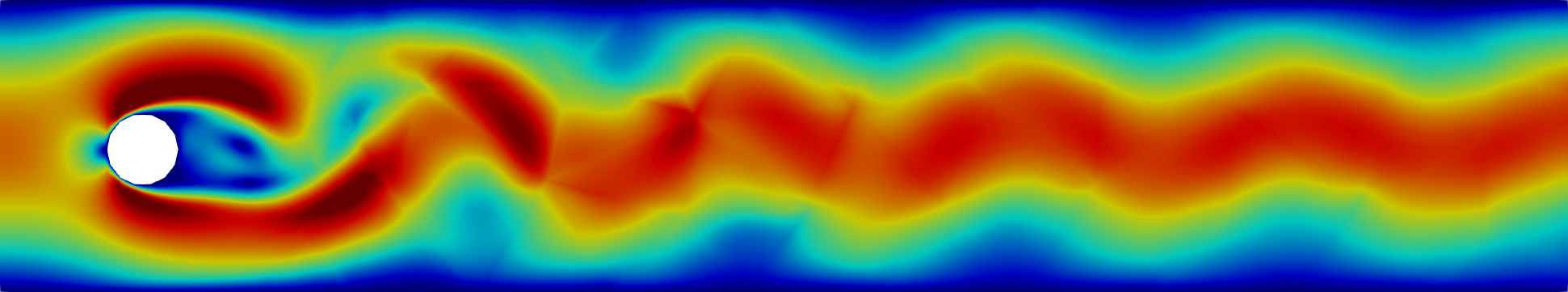}
  \end{center}
    \caption{\corrections{Flow around a cylinder showing magnitude of velocity.
    On the top left we show \changes{the results from using a} second-order \changes{space for the} velocity and piecewise constant pressure on a coarse grid with 260 triangles which itself is shown on the top right figure.
    On the middle right figure we show the \changes{the results from} using a fourth-order velocity space and piecewise constant pressure on the same coarse grid.
    The middle left figure shows results with the same space as top left but on a grid with 1410 triangles having the same resolution around the cylinder but a higher resolution downstream.
    Bottom left figure shows results using a second-order Taylor Hood space on the same coarse grid shown in the top right figure and finally, the bottom right figure shows results using a fourth-order Taylor Hood space on the same coarse grid.}}
\end{figure}

We use the divergence free space to solve the Navier-Stokes equations
for an incompressible flow with $\nu=0.001$\corrections{: we seek $(u(t),p(t))$ such that:} 
\begin{align*}
  \partial_t u + u\cdot\nabla u - \nu\Delta u + \nabla p &= 0~,
  \quad {\rm div} u = 0 \quad \text{in } \Omega.&
\end{align*}
Here, $p$ is the pressure and for \corrections{the velocity} $u$ we prescribe the usual initial and boundary conditions for a flow around a cylinder with radius $0.05$ located at $(0.2,0.2)$ in
the domain $[0,2.2]\times[0,0.41]$ \corrections{see e.g. \cite{john2002higher,schafer1996benchmark}}. 
To discretise this system we use the divergence free space for the velocity and a piecewise constant approximation for the pressure, which is the compatible space containing the divergence of the discrete velocity space. 
We use a simple discretisation in time based on a semi-implicit method
\begin{align*}
  \frac{1}{\tau}u^{n+1} - \nu\Delta u^{n+1} + \nabla p^{n+1}
        &= \frac{1}{\tau}u^n - u^n\cdot\nabla u^n~,\quad {\rm div} u^{n+1} = 0
\end{align*}
with a time step $\tau = 6.25 e-4$. The resulting saddle point problem is
solved using an Uzawa-type algorithm in each time step. 
We solve on a triangular grid using different polynomial
degrees and compare with results using a more standard Taylor-Hood space.
Results are shown in \cref{fig: navier stokes}.

\subsection{Willmore flow of graphs}
In our final example we study the minimisation of the Willmore energy of a surface,
in the case where the surface is given by a graph over a flat domain $\Omega$.
The corresponding Euler-Lagrange equations can be rewritten as a fourth-order problem for a function $u$ defined over $\Omega$. 
We use a second-order two stage implicit Runge-Kutta method as suggested in \cite{deckelnick_c1finite_2015}. 
The resulting problem is a system of two nonlinear fourth-order partial differential equations for the two Runge-Kutta stages.

\begin{figure}[!ht]\label{fig: willmore}
\centering
\includegraphics[width=0.75\textwidth]{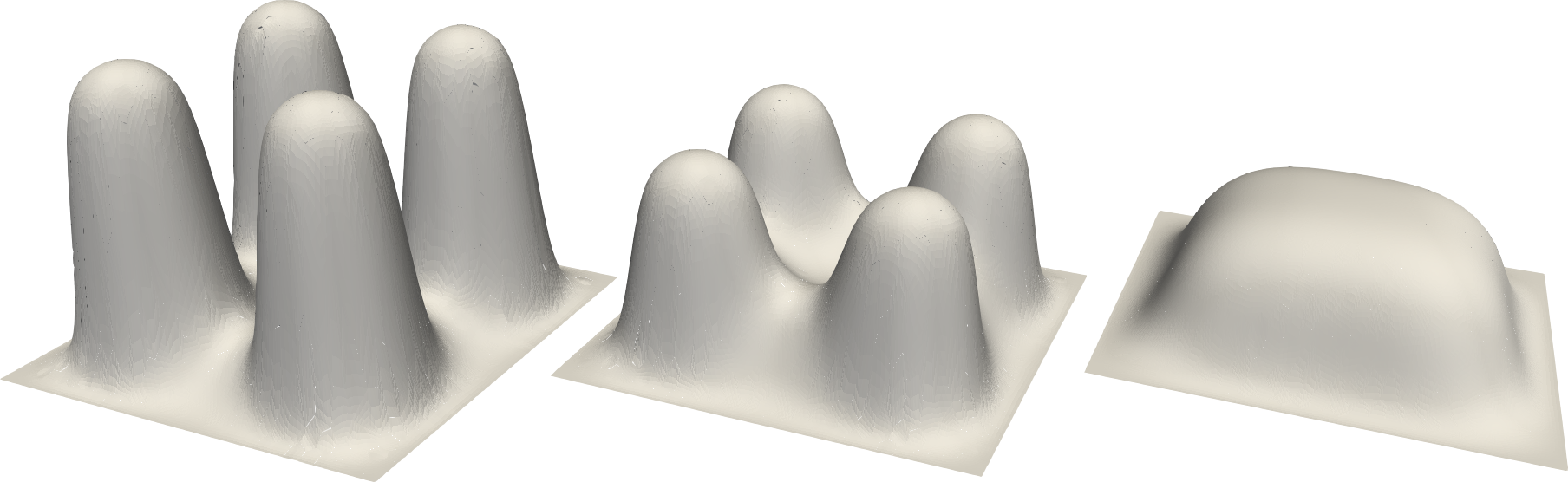} 
\caption{Evolution of a graph under Willmore flow at times
$t=1e-4,4e-4,7e-4$ (from left to right) using the $H^2$-conforming space of order $\polOrder=4$. Results with the nonconforming space (also considered in \cite{10.1093/imanum/drab003}) of the same order are indistinguishable.}
\end{figure}

As detailed for example in \cite{deckelnick_c1finite_2015}, the Willmore functional for the graph of a function ${u \in W^{2,\infty}(\Omega)}$ is given by 
\begin{align*}
  W(u) = \frac{1}{2}\int_{\Omega} [ \ E(\nabla u) : \nabla^2  u \ ]^2 \, \dx~,\;\text{with}\;
  E_{ij}(w) := \frac{1}{(1+|w|^2)^{\frac{1}{4}}} \big( \delta_{ij} - \frac{w_i w_j}{1+|w|^2} \big)
\end{align*}
for $i,j=1,2$, and $w \in \R^2$.
\corrections{We initialise the gradient descent algorithm with $u(x,y)=\big(\sin(2\pi x)\sin(2\pi y)\big)^2$} and we use the time step $\tau = 5e-6$ with the $H^2$-conforming space from \Cref{sec: H2 conforming vem}.
\cref{fig: willmore} shows the evolution of the graph
on a Voronoi grid with $800$ cells.

\section{Conclusion}\label{sec: conclusion}
In this paper, we have presented a framework for implementing general virtual element spaces in two space dimensions \corrections{and discussed what we believe is a straightforward extension of the framework concepts to three dimensions.}  
As is usual with VEM schemes, the definition of \emph{projection} operators is crucial \corrections{in order} to setup the discrete bilinear forms. In our approach the definition of the gradient and hessian projections are independent of the space \corrections{and are} based on value projections on the elements and skeleton of the grid.
We introduced a VEM tuple for encapsulating all the necessary building blocks for computing these projection operators, a concept which aimed to mimic the FEM triple from the finite element setting.
These building blocks included basis sets, dof sets, and constraint sets needed to construct the value projection operators on the elements and skeleton. 
\corrections{Additionally,} basis sets can be provided for the gradient and hessian projection.
Our starting point for the construction of the value projections is \refereecommentsvisible{constrained} least squares problems.
Projections for the higher order derivatives are then defined independently of the space.
With examples we showed how to construct different VEM spaces with additional properties such as $H^k$-conforming spaces for $k=1,2$, divergence free as well as curl free spaces.
Our approach has the added benefit of encapsulating extra properties of the space through the value projection.
\changes{This avoids requiring further projection operators for e.g. the divergence, thus simplifying integration into existing frameworks.}

\changes{One major advantage of our framework is that it can be easily integrated into an existing finite element package.
As already mentioned the main advantage of the presented formulation is that no special projections depending on the underlying PDE are utilised thus minimising the changes required to existing software frameworks.}
We demonstrated this in two space dimensions within the \textsc{Dune} software framework and presented a handful of numerical experiments to showcase the wide variety of VEM spaces which can be utilised.
To the best of our knowledge, this is the first available implementation to include such a vast collection of VEM spaces including but not limited to, spaces for fourth-order problems and nonlinear problems.
The software is free and open source and the $\texttt{dune-vem}$ module can be easily installed using PyPI.

\bibliographystyle{siamplain} 
\bibliography{compatible_VEM}

\end{document}